\documentclass[reqno,11pt]{article}
\usepackage{graphicx, amsmath, amstexnb, amsthm}
\usepackage{amssymb}  

\def\captionspace{}

\input{stochbif.sty}

\begin{document}


\title{On the noise-induced passage\\through an unstable periodic orbit I:\\ 
Two-level model}
\author{Nils Berglund and Barbara Gentz}
\date{}   

\maketitle

\begin{abstract}
\noindent
We consider the problem of stochastic exit from a planar domain, whose
boundary is an unstable periodic orbit, and which contains a stable
periodic orbit. This problem arises when investigating the distribution of
noise-induced phase slips between synchronized oscillators, or when
studying stochastic resonance far from the adiabatic limit. We introduce a
simple, piecewise linear model equation, for which the distribution of
first-passage times can be precisely computed. In particular, we obtain a
quantitative description of the phenomenon of cycling: The distribution of
first-passage times rotates around the unstable orbit, periodically in the
logarithm of the noise intensity, and thus does not converge in the
zero-noise limit. We compute explicitly the cycling profile, which is
universal in the sense that in depends only on the product of the period of
the unstable orbit with its Lyapunov exponent. 
\end{abstract}

\leftline{\small{\it Date.\/} August 19, 2003.}
\leftline{\small 2000 {\it Mathematical Subject Classification.\/} 
37H20, 60H10 (primary), 34F05 (secondary)}
\noindent{\small{\it Keywords and phrases.\/}
Stochastic exit problem, diffusion exit, first-exit time, large
deviations, metastability, level-crossing problem, limit cycle,
synchronization, phase slip, cycling, stochastic resonance.}  


\section{Introduction}
\label{sec_in}


One of the remarkable effects of additive noise on a deterministic
dynamical system is to enable transitions between otherwise isolated
attractors. This phenomenon, known as {\em activation\/} in physics and
chemistry, is at the origin of the so-called {\em stochastic exit problem}.
Consider a stochastic differential equation (SDE) of the form 
\begin{equation}
\label{in1}
\6x_t = f(x_t)\6t + \sigma\6W_t.
\end{equation}
Given an attractor $\cA$ of the deterministic system $\dot x=f(x)$ and a
domain $\cD$ containing $\cA$ (usually $\cD$ is taken to be positively
invariant under the deterministic flow), the exit problem consists in
characterizing
\begin{itemiz}
\item	the distribution of the random time $\tau =
\inf\setsuch{t>0}{x_t\notin\cD}$ at which paths of the SDE, starting in some
initial point in $\cD$, leave the domain;
\item	the distribution of the exit location $x_\tau$ on the boundary
$\partial\cD$ of the domain. 
\end{itemiz}
The best understood situation is the one where $f=-\nabla U$ derives from a
potential $U$, and the attractor $\cA$ is simply the bottom $x^\star$ of a
potential well. Assume that $\cD$ contains no other equilibrium point in its
interior, and a (non-degenerate) saddle point $x_1$ of $U$ on its
boundary. The expectation of the first-exit time $\tau$ is then given by 
\begin{equation}
\label{in2}
\expec{\tau} = c(\sigma) \e^{[U(x_1)-U(x^\star)]/2\sigma^2}. 
\end{equation}
The form of the exponent was already known to Arrhenius~\cite{Arrhenius},
while the value of the prefactor $c(\sigma)$, which depends on the curvature
of $U$ at $x^\star$ and $x_1$, was determined by Eyring and Kramers in the
one-dimensional case~\cite{Eyring,Kramers}. 

Many refinements of this result exist. On the one hand, the theory of large
deviations allows to compute the exponent in a much more general setting,
with a drift term not necessarily deriving from a potential~\cite{FW}. In
this case, $U$ is replaced by the notion of a {\em quasipotential}. On
the other hand, the computation of the prefactor $c(\sigma)$ in $n$
dimensions has been addressed using methods of perturbation theory
(see for instance \cite{Azencott,FJ}). Precise results on the
relation between the expected first-exit time, capacities and
the small eigenvalues of the generator of the diffusion have been
obtained recently for drift coefficients deriving from a
potential~\cite{BEGK,BGK}.   

The distribution of $\tau$ approaches an exponential one in the small-noise
limit~\cite{Day1,BGK}, and the exit location $x_\tau$ is strongly
concentrated near the saddle $x_1$~\cite{FW}. Another related property
is that for a 
fixed time $t\ll\expec{\tau}$, the probability to leave $\cD$ before
time~$t$ is exponentially small, of the order $t/\expec{\tau}$. This
phenomenon is known as {\em metastability}, because the state $x_t$ may
spend extremely long time spans near local minima of the potential, while
the stationary density has most of its mass concentrated near the global
potential minima. 

A more difficult situation arises, in two-dimensional systems, when the
boundary of $\cD$ is an unstable periodic orbit of the deterministic flow.
In that case, large deviations theory cannot be applied directly, because
$\partial\cD$ is a so-called characteristic boundary. 
The situation where $\partial\cD$ is an unstable periodic orbit and $\cA$
an equilibrium point has been considered by Day~\cite{Day5,Day3} and by
Maier and Stein~\cite{MS4}. Day~\cite{Day3} proved that the distribution of the
first-exit location $x_\tau$ displays a striking behaviour, called {\em
cycling\/}: As the noise intensity $\sigma$ decreases, the density of
$x_\tau$ rotates around the boundary, as a function of $\abs{\log\sigma}$.
Using the concept of the most probable exit path (MPEP) and WKB
approximations, Maier and Stein also found a periodic dependence on
$\abs{\log\sigma}$ for the rate of escape per unit time through
$\partial\cD$. The intuition behind their estimates is that the MPEP
spirals geometrically towards $\partial\cD$, until it reaches a
neighbourhood of order $\sigma$ of $\partial\cD$ and escape becomes
likely. Where and when this happens depends periodically on 
$\abs{\log\sigma}$. 

A related situation, which has not yet been analysed in detail, arises when
the boundary of $\cD$ is an unstable periodic orbit but $\cA$ is a stable
periodic orbit. This situation is important in applications:

\begin{enum}
\item	When studying the dynamics of two coupled, slightly different phase
oscillators, the onset of {\em synchronization\/} is known to correspond to
a saddle--node bifurcation of periodic orbits.  Below a threshold coupling
strength, the motion is typically quasiperiodic, meaning that the
oscillators are not synchronized. Above this threshold, a pair of periodic
orbits of opposite stability appears, the stable one corresponding to a
synchronized state in which the phase difference between the oscillators
may vary periodically but is bounded (see for instance~\cite{PRK}).

Adding noise to the system will cause {\em phase slips\/} to occur, in which
the phase difference changes by $2\pi$. Such a phase slip
necessarily involves crossing the unstable periodic orbit, so that the
determination of the distribution of phase slips requires the determination
of the distribution of first-passage times at an unstable orbit. 

\item	The phenomenon of {\em stochastic resonance\/} occurs for instance
when a double-well potential is forced periodically in
time~\cite{BSV}. Noise-induced transitions between potential
wells are more favourable from the shallower to the deeper well, when
the barrier between them is lowest. As a result, typical paths of the
system will contain a periodic component, see for instance~\cite{Fox,WM,GHM}.

On the mathematically rigorous level, the situation is relatively well
understood in the adiabatic case, that is for slow forcing, when the paths
spend most of the time near the bottom of a potential
well~\cite{Freidlin1,BG2}. When the forcing is not slow, however, it is
easy to see that the deterministic system still admits three periodic
solutions, two stable ones oscillating around the potential wells, and an
unstable one oscillating around the saddle. Hence the investigation of
transition between wells again involves the understanding of first-passage
times at the unstable orbit.  
\end{enum}

One may point out that in the case of synchronization, the phase space has
the topology of a torus, while for stochastic resonance, it has the
topology of a cylinder. However, both situations have in common the fact
that a path, starting near a stable periodic orbit, has to cross an
unstable periodic orbit for a transition to become possible. 

In the present work, we focus on the dynamics of the paths up to their first
crossing of $\partial\cD$. Our aim is to give a precise characterization of
the distribution of exit times and locations, in particular in the
metastable regime. In order to highlight the mechanism responsible for
cycling, we concentrate here on a simplified, piecewise linear model
equation, which can be solved exactly to leading order. The general case
will be discussed in a forthcoming publication~\cite{BG8}. 

In our model, we neglect the diffusion in the longitudinal direction, so
that the exit location is determined by the exit time, modulo the period of
the orbit. The main result, Theorem~\ref{thm_mrr}, gives an explicit
expression for the density of the first-exit time. If $T$ is the period
of the unstable orbit, and $\lambda$ its Lyapunov exponent, for a large
range of metastable times $t$ such that $2\abs{\log\sigma}\ll \lambda t\ll
\e^{\const/\sigma^2}$, this density is given by 
\begin{equation}
\label{in3}
p_+(t) \simeq \const \sigma \theta'(t) 
P\biggpar{\frac{\abs{\log\sigma}-\theta(t)}{\lambda T}} \e^{-R^2/2\sigma^2}. 
\end{equation} 
Here $R^2$ describes the exponential rate of escape provided by large
deviations theory. The {\em cycling profile\/} $P(x)$ is an explicitly known,
universal function, depending only on $\lambda T$. The model-dependent
intrinsic time $\theta(t)$, which satisfies $\theta(t+T)=\theta(t)+\lambda
T$, describes the \lq\lq velocity\rq\rq\ $1/\theta'(t)$ with which the
cycling profile rotates around the orbit. 

The regime $\lambda t\leqs2\abs{\log\sigma}$ is transient, in the sense
that paths have not yet reached their typical spreading, and thus $p_+(t)$
is smaller than~\eqref{in3} by a factor $\e^{-\const \e^{-\lambda
t}/\sigma^2}$. 
This initial phase is not observed in~\cite{MS4}, where the authors
artificially create a stationary regime by reinjecting escaped paths into
the attractor. 

Though we do not treat in detail the asymptotic regime
$t\gg\e^{\const/\sigma^2}$, results by Day~\cite{Day3} imply that the
superposition 
$\sum_{k\geqs0}p_+(t+kT)$, which does not take into account the winding
number~$k$ of paths around the unstable orbit, has a similar behaviour
as~\eqref{in3}. In fact, we expect~\eqref{in3} to hold for times larger than
$\e^{\const/\sigma^2}$, with an additional factor, slowly decaying like 
$\exp\set{-t\e^{-R^2/2\sigma^2}}$.

\subsection*{Acknowledgements}

The present work was stimulated by discussions with Arkady Pikovsky.  We
thank Peter M\"orters for bringing recent results on first-passage time
densities to our attention.

N.B. thanks the WIAS for kind hospitality and financial support.  B.G.
thanks the Universit\'e de Toulon and the CPT--CNRS Luminy for hospitality.
Financial support by the ESF Programme {\it Phase Transitions and
Fluctuation Phenomena for Random Dynamics in Spatially Extended Systems
(RDSES)\/} is gratefully acknowledged.


\section{Model and results}
\label{sec_mr}


\subsection{Periodic orbits and coordinate systems}
\label{ssec_mrdc}

Consider a two-dimensional ordinary differential equation $\dot x=f(x)$
admitting a stable periodic orbit enclosed by an unstable one. We assume
that the domain $\cS$ lying between the orbits has the topology of an
annulus and contains no invariant sets. In that case, one can choose polar-like
coordinates $(r,\ph)$, such that $\ph$ is $2\pi$-periodic and $\dot\ph>0$ in
a neighbourhood of $\cS$. 

In the deterministic case, it is customary to use $\ph$ as new independent
variable, to obtain a one-dimensional non-autonomous system
\begin{equation}
\label{mrdc1}
\dtot r\ph = f_r(r,\ph).
\end{equation}
Furthermore, it is possible to choose $r$ in such a way that the stable
periodic orbit corresponds to $r=-1$, and the unstable one to $r=+1$. 

\begin{example}
\label{ex_periodic1}
{\bf Synchronization} \hfill

\noindent
The onset of synchronization between two weakly coupled phase oscillators
can be viewed as a saddle--node bifurcation of periodic orbits. A 
normal-form analysis shows that close to the bifurcation point, the
dynamics is governed to leading order by an equation of the form 
\begin{equation}
\label{mrdc2}
\dtot r\ph = r^2 - \eps c(\ph),
\end{equation}
where $\eps$ is the bifurcation parameter. If $c(\ph)>0$ for all $\ph$, a
straightforward analysis of the Poincar\'e map shows the existence, for
$\eps>0$, of two periodic orbits $r=r_+(\ph)$ and $r=r_-(\ph)$ of opposite
stability, separated by a distance of order $\sqrt\eps$. The linear
transformation 
\begin{equation}
\label{mrdc3}
r = \frac{r_+(\ph)+r_-(\ph)}2 + \frac{r_+(\ph)-r_-(\ph)}2 y
\end{equation}
yields the equation 
\begin{equation}
\label{mrdc4}
\dtot y\ph = \frac12\bigbrak{r_+(\ph)-r_-(\ph)} \bigpar{y^2-1}. 
\end{equation}
The periodic orbits are now located in $y=\pm1$. Note that the stable
(unstable) orbit is  attracting (repelling)
more strongly for those values of $\ph$ for which the orbits in the
original system are further apart. 
\end{example} 

When noise is added to the system, it will in general affect both the
transversal $r$- and the longitudinal $\ph$-direction, so that one is led to
analyse the system 
\begin{equation}
\label{mrdc5}
\begin{split}
\6r_t &= f_r(r_t,\ph_t) \6t + \sigma g_r(r_t,\ph_t) \6W_t, \\
\6\ph_t &= f_\ph(r_t,\ph_t) \6t + \sigma g_\ph(r_t,\ph_t) \6W_t,
\end{split}
\end{equation}
possibly with $f_\ph=1$. Here we consider $\sigma$ as a small parameter
controlling the noise intensity, while $g_r$ and $g_\ph$ are fixed functions
(of order one). 

\begin{example}
\label{ex_periodic2}
{\bf Stochastic resonance} \hfill

\noindent
A classical example of a system showing stochastic resonance consists of an
overdamped particle in a periodically forced double-well potential,
perturbed by additive noise. For a Ginzburg-Landau potential, the equation
reads 
\begin{equation}
\label{mrdc7}
\6r_t = \bigbrak{r_t - r_t^3 + A\cos(2\pi t/T)}\6t + \sigma\6W_t.
\end{equation}
Again, using a Poincar\'e section, one easily proves the existence of one
unstable and two stable periodic orbits (for any value of the period $T$).
If $r_-(t)$ denotes, say, the lower stable orbit, and $r_+(t)$ the unstable
one, a transformation of the form~\eqref{mrdc3} yields
Equation~\eqref{mrdc5} with $f_\ph=1$, $g_\ph=0$ (that is, $\ph=t$), and
$g_r$ depending only on $\ph$. 
\end{example}

The general case of Equation~\eqref{mrdc5} will be discussed in a forthcoming
work \cite{BG8}. Here we shall concentrate on a simplified model which
focuses on the main mechanism responsible for the oscillatory behaviour of
the first-exit time. 


\subsection{Simplified two-level model}
\label{ssec_mrs}

The main motivation to introduce the two-level model is the fact that the
distribution of first-exit times will be mainly determined by the
dynamics near the unstable periodic orbit. The dynamics in the
remaining phase space can thus be modelled by the simplest possible
equation, that is, a linear one. More precisely, we will
simplify~\eqref{mrdc5} by  
\begin{itemiz}
\item	neglecting the term $g_\ph$, whose effect is a slow
diffusion on a time scale $1/\sigma^2$;
\item	neglecting the $r$-dependence of $g_r$, which is not important near
the unstable orbit;
\item	replacing $f_r$ by a piecewise linear function of $r$. 
\end{itemiz}
With these approximations, we arrive at the system
\begin{equation}
\label{mrdc6}
\6y_t = f(y_t,t) \6t + \sigma g(t) \6W_t, 
\end{equation}
where $f$ and $g$ are periodic in $t$, and we may assume that $f(\pm1,t)=0$
for all $t$. Here time $t$ is identified with the angle $\ph$ lifted to the
real axis. 

\begin{figure}
\centerline{\includegraphics*[clip=true,width=125mm]{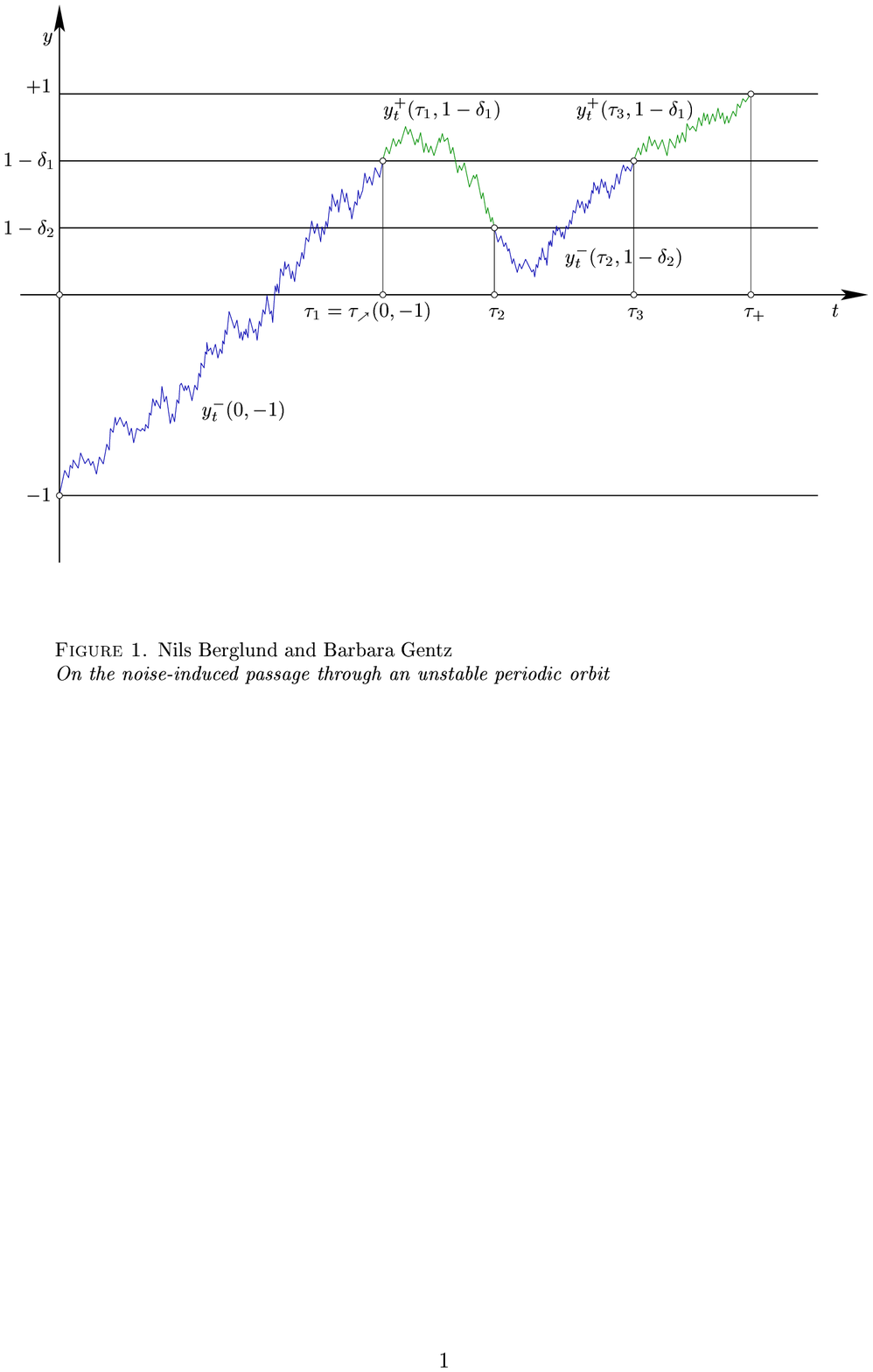}}
 \captionspace
 \caption[]
 {The process $y_t$ is defined by switching between two linear
 processes $y^-_t$ and $y^+_t$ each time either the level $1-\deltap$ is
 crossed from below or the level $1-\deltam$ is crossed from above.
 }
\label{fig1}
\end{figure}
In order to avoid certain technical difficulties when dealing with the
stochastic process, we will actually switch between two linear equations
defined in slightly overlapping regions. These equations are 
\begin{equation}
\label{mrs1}
\begin{split}
\6 y^-_t &= - a(t) (y^-_t+1) \6t + \sigma g(t) \6W_t, \\
\6 y^+_t &= \phantom- a(t) (y^+_t-1) \6t + \sigma g(t) \6W_t,
\end{split}
\end{equation}
where $a(t)$ and $g(t)$ are $T$-periodic, positive functions, which are
bounded away from zero (detailed assumptions will be given in
Section~\ref{ssec_mrr} below). We denote by $y^-_t(t_0,y_0)$ and 
$y^+_t(t_0,y_0)$ the solutions of these equations with initial conditions
$y^-_{t_0}=y_0$ or  $y^+_{t_0}=y_0$, respectively.  The stable orbit
located at $y=-1$ and the unstable orbit at $y=+1$ have Lyapunov exponents
$\mp\lambda$, where 
\begin{equation}
\label{mrs2}
\lambda = \frac{\alpha(T)}T, 
\qquad
\text{with}
\qquad
\alpha(t) = \int_0^t a(s)\6s.
\end{equation}
The switching between the processes occurs upon reaching levels $1-\deltap
\in(0,1)$ from below and $1-\deltam\in(0,1-\deltap)$ from above
(see~\figref{fig1}). More precisely, consider the stopping times  
\begin{equation}
\label{mrs3}
\begin{split}
\tauup &= \tauup(t_0,y_0) 
= \inf\setsuch{t>t_0}{y^-_t(t_0,y_0)>1-\deltap} \in[t_0,\infty], \\
\taudown &= \taudown(t_0,y_0) 
= \inf\setsuch{t>t_0}{y^+_t(t_0,y_0)<1-\deltam} \in[t_0,\infty]. 
\end{split}
\end{equation}
Then the process $\set{y_t}_{t\geqs0}$ is defined in the following way:
\begin{equation}
\label{mrs4}
y_t = 
\begin{cases}
y^-_t(0,-1) 
&\text{for $\,\phantom{\tau_1}\llap{0}\leqs t\leqs \tau_1 = \tauup(0,-1)$}, \\
y^+_t(\tau_1,1-\deltap) 
&\text{for $\,\tau_1\leqs t\leqs \tau_2 = \taudown(\tau_1,1-\deltap)$}, \\
y^-_t(\tau_2,1-\deltam) 
&\text{for $\,\tau_2\leqs t\leqs \tau_3 = \tauup(\tau_2,1-\deltam)$},
\end{cases}
\end{equation}
and so on. We are interested in determining the distribution of the
first-passage time at the unstable orbit, namely in the
distribution of 
\begin{equation}
\label{mrs5}
\tau_+ = \inf\setsuch{t>0}{y_t>+1}. 
\end{equation}
We use the notation 
\begin{equation}
\label{mrs6}
p_+(t) = \dpar{}{t} \probin{0,-1}{\tau_+\leqs t}
\end{equation}
for the density of $\tau_+$, where the superscript in
$\probin{0,-1}{\cdot}$ indicates that we start at time $t=0$ on the
stable orbit in $y=-1$. For brevity, we shall call $p_+(t)$ the
first-passage density of $y_t$ to $+1$.  

Note that since $\ph$ is proportional to $t\pmod{T}$, the first-exit time
$\tau_+$ directly allows to determine the first-exit location
$\ph_{\tau_+}$. However, $\tau_+$ actually contains more information since
it also keeps track of the number of revolutions (or winding number) of the
path around the stable orbit. 

In the physics literature, one often considers the rate of escape per unit
time, defined as minus the time-derivative of the probability to be inside
the domain $\cD$. The difference between rate of escape and first-passage
density is that the former counts negatively the paths which have left $\cD$
before time $t$, but returned into $\cD$ by time $t$, while these paths are
not counted by the first-passage density. If the system is symmetric with
respect to the boundary $\cD$, the reflection principles implies that the
rate of escape is equal to half the first-passage density.  


\subsection{Main results}
\label{ssec_mrr}

\begin{figure}
\centerline{\includegraphics*[clip=true,width=100mm]{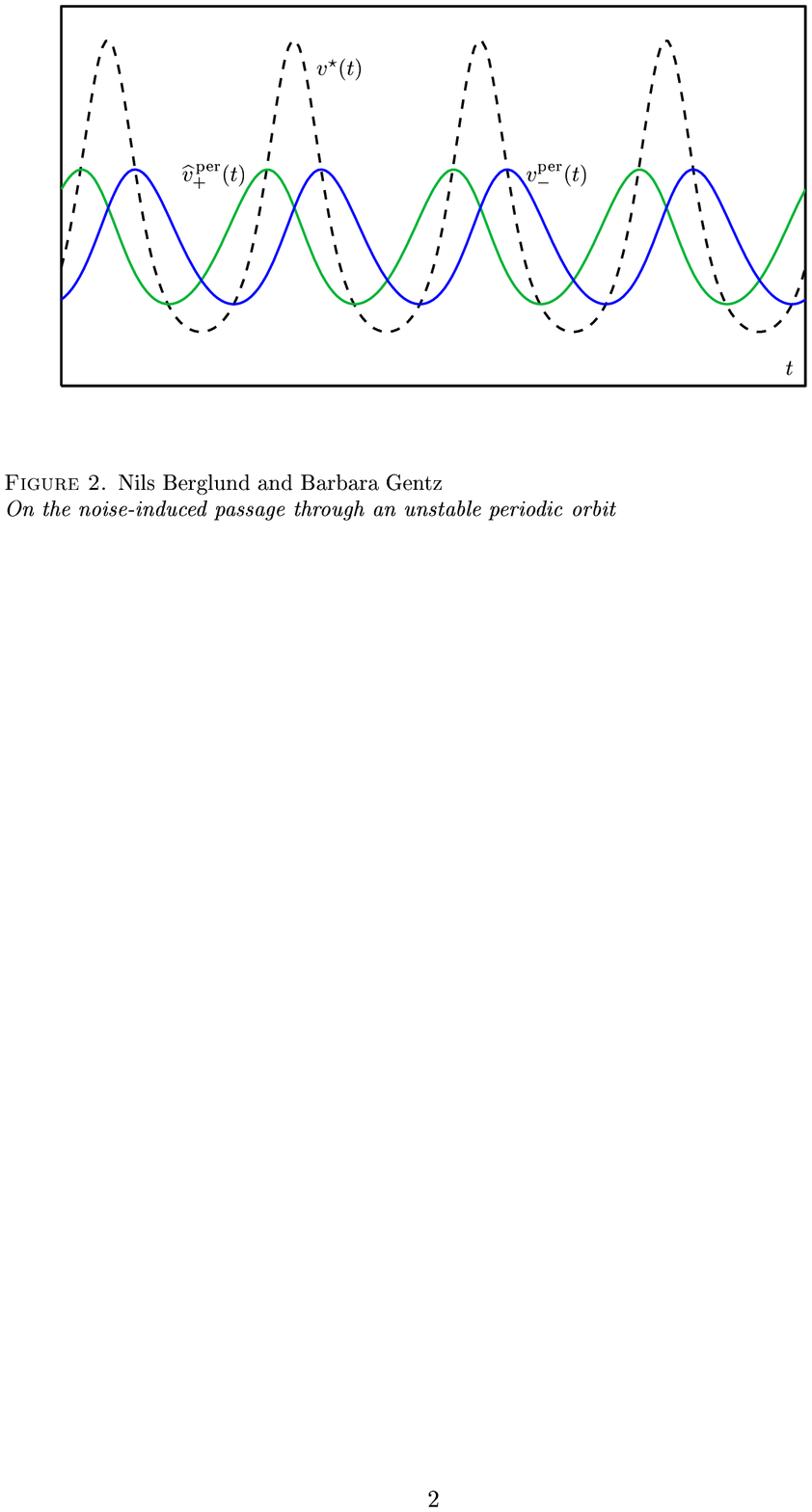}}
 \captionspace
 \caption[]
 {The periodic functions $v^\star(t)$, $\vper_-(t)$ and $\vhatper_+(t)$.}
\label{fig2}
\end{figure}

Four periodic functions will play an important r\^ole in the statement of
the results. They are given by 
\begin{align}
\label{mrr1a}
\vper_-(t) &= \frac1{1-\e^{-2\lambda T}} \int_t^{t+T} \e^{-2\alpha(t+T,s)}
g(s)^2\6s, \\
\label{mrr1b}
\vhatper_+(t) &= \frac1{\e^{2\lambda T}-1} \int_t^{t+T} \e^{2\alpha(t+T,s)}
g(s)^2\6s, \\
\label{mrr1c}
\rhoper(t)^2 &= \frac{\deltap^2}{\vhatper_+(t)} 
+ \frac{(2-\deltap)^2}{\vper_-(t)}, \\ 
\label{mrr1d}
v^\star(t) &= \frac{g(t)^2}{2a(t)}, 
\end{align}
where $\alpha(t,s)=\alpha(t)-\alpha(s)$. The first two functions are
directly related to the variances of $y^-_t(t_0,y_0)$ and $y^+_t(t_0,y_0)$.
These are independent of the initial condition $y_0$, and have respective
values
\begin{align}
\label{mrr2a}
\sigma^2 v_-(t,t_0) &= \sigma^2 \bigbrak{\vper_-(t) -
\e^{-2\alpha(t,t_0)}\vper_-(t_0)}, \\
\label{mrr2b}
\sigma^2 v_+(t,t_0) &= \sigma^2  
\bigbrak{\e^{2\alpha(t,t_0)}\vhatper_+(t_0) - \vhatper_+(t)}.
\end{align}
The function $v^\star(t)$ allows to determine the qualitative behaviour of
$\vper_-(t)$ and $\vhatper_+(t)$. Indeed, since $\vper_-(t)$ satisfies the
differential equation $\dot v = -2a(t)v + g(t)^2$, it is increasing whenever
it lies below $v^\star(t)$ and decreasing whenever it lies above.
Similarly, $\vhatper_+(t)$ is decreasing when it lies below $v^\star(t)$
and increasing when it lies above (\figref{fig2}). This shows in particular
that $\vper_-(t)$ and $\vhatper_+(t)$ always lie between the maximum $\bar
v$ and the minimum $\underline v$ of $v^\star(t)$, and that they reach
their extremal values when crossing $v^\star(t)$. One can also see from the
graphical representation that smaller values of the period $T$ lead to
smaller amplitudes of $\vper_-(t)$ and $\vhatper_+(t)$. 

Our main assumptions are the following:

\bigskip\noindent
{\bf Hypotheses.}
\label{basic_assump}
\begin{enumH}
\item[H1.]	The functions $a(t)$ and $g(t)$ are twice continuously
differentiable, positive, $T$-periodic and bounded away from zero.
\item[H2.]	$v^\star(t)$ has exactly one maximum and one minimum in $[0,T)$,
with values $\bv>\underline v>0$. 
\item[H3.]	There is a constant $\Delta>0$ such that $\vper_-(t)\leqs
2v^\star(t)(1-\Delta)$ and $\vhatper_+(t)\leqs 2v^\star(t)(1-\Delta)$ for
all $t$. 
\item[H4.]	$\deltam/(2-\deltam) \leqs \sqrt{\underline v/\bv}$. 
\item[H5.]	$\rhoper(t)$ has exactly one minimum in $[0,T)$ at a time
$t=s^\star$, which is quadratic. We set 
\begin{equation}
\label{mrr3}
R = \rhoper(s^\star) = \inf_{t\in[0,T)} \rhoper(t). 
\end{equation}
\end{enumH}
\smallskip

Note that Hypothesis~H3 amounts to requiring that the most probable exit paths
of the processes $y^-_t$ and $y^+_t$ can cross at most once the levels
$1-\delta_1$ and $1-\delta_2$ (see the Appendix). It is always satisfied
when $T$ is large enough, because then the variances $\vper_-$ and
$\vhatper_+$ track $v^\star$ adiabatically. Is is also satisfied if
$\bv<2\underline v(1-\Delta)$. Hypothesis~H4 ensures that $y^-_t$ is
dominated by $y^+_t$ near the unstable orbit. We make Hypothesis~H5 mainly
in order to simplify the presentation, it is in fact sufficient to require
that the deepest minimum of $\rhoper(t)$ be quadratic, which is generically
true if $v^\star(t)$ is nontrivial. 

Under Assumptions H1--H5, the following result on the first-passage
density $p_+(t)$ holds whenever $\sigma$ is small enough.

\begin{theorem}
\label{thm_mrr}
There exists a $\sigma_0>0$ such that for all $\sigma\leqs\sigma_0$, 
the density $p_+(t)$ of $\tau_+$ is given by 
\begin{equation}
\label{mrr4}
p_+(t) = c(t,\sigma) \e^{-R^2/2\sigma^2}, 
\end{equation}
where $R$ is defined in~\eqref{mrr3} and the prefactor $c(t,\sigma)$ depends
on $t$ and $\sigma$ in the following way: 
\begin{enum}
\item	{\bf Transient regime:} If\/ $0\leqs \alpha(t) <
2\abs{\log\sigma}$, then 
\begin{equation}
\label{mrr5}
c(t,\sigma) \leqs \const
\frac1{\sigma^2} \exp\biggset{-\frac{L \e^{-\alpha(t)}}{\sigma^2}},
\end{equation}
where $L$ is a positive constant. 

\item	{\bf Metastable regime:} Let $\Delta_0 = \Delta/(1-\Delta)
\wedge (\delta_2-\delta_1)/\delta_1\wedge1$.\footnote{We write
$a\wedge b$ to denote the minimum of two real numbers $a$ and $b$.}\ %
There exists a constant $\beta>0$
such that for $t$ satisfying\/ $2\abs{\log\sigma} \leqs \alpha(t) \ll
\sigma^3\e^{\beta\Delta_0^2/2\sigma^2}$, 
\begin{equation}
\label{mrr6}
c(t,\sigma) = 
\sigma C_0(s^\star) \theta'(t)
P\biggpar{\frac{\abs{\log\sigma}-\theta(t)}{\lambda T}} 
\bigbrak{1 + r(\sigma)},
\end{equation}
where $P(x)>0$ is an explicitly known periodic function with period $1$, see
Equations~\eqref{mrr8} and~\eqref{mrr8b} below, and
\begin{align}
\label{mrr6b}
\theta(t) &= \alpha(t,s^\star) - \frac12\log\vhatper_+(t)
- \log\frac{2-\deltap}{\vhatper_+(s^\star)},\\ 
\label{mrr6c}
\theta'(t) &= \frac12 \frac{g(t)^2}{\vhatper_+(t)}. 
\end{align}
$C_0(s^\star)$ is a constant given by 
\begin{equation}
\label{mrr7}
C_0(s^\star) = 4\,\frac{2-\deltap}{\deltap}
\frac{g(s^\star)^2}{\sqrt{\pi\sdpar{}{ss}(\rhoper(s^\star)^2)}}
\frac{\vhatper_+(s^\star)^{1/2}}{\vper_-(s^\star)^{3/2}}
\biggbrak{1-\frac{\vper_-(s^\star)}{2v^\star(s^\star)}},
\end{equation}
and the error term satisfies
\begin{equation}
\label{mrr7.5}
r(\sigma) = \BigOrder{\sigma + \frac{1}{\sigma^2} \e^{-\alpha(t)}}.
\end{equation}

\item	{\bf Asymptotic regime:} If\/ $\alpha(t) \geqs
\const \e^{R/2\sigma^2}$, most paths will have crossed the
unstable orbit as least once, and thus the density decays. 
\end{enum}
\end{theorem}

The periodic function $P(x)$ appearing in~\eqref{mrr6} has the following
expressions: 
\begin{equation}
\label{mrr8}
P(x) = \sum_{\l=-\infty}^{\infty} A(\lambda T(\l-x)), 
\qquad
\text{with}
\quad
 A(x) = \frac12 \e^{-2x}\exp\Bigset{-\frac12 \e^{-2x}},
\end{equation}
which is particularly useful for large $T$, and the Fourier series 
\begin{equation}
\label{mrr8b}
P(x) = 
\sum_{q=-\infty}^{\infty} 
\Phat(q)\e^{2\pi\icx qx},
\qquad
\text{with}
\quad 
\Phat(q) = \frac{1}{2\lambda T} \frac1{2^{\pi\icx q/\lambda T}}
\Gamma\biggpar{1 - \frac{\pi\icx q}{\lambda T}},
\end{equation}
where $\Gamma$ is the Euler Gamma function. This series converges quickly
when $T$ is small. 
 
Before discussing the implications of this result, let us briefly
sketch the proof, the details of which are given in Sections~\ref{sec_re}
to~\ref{sec_xl}. 

A first step is to determine the density $\psi_-(s,0)$ of the first-passage
time at $1-\deltap$, when the first switching occurs. We will show in
Proposition~\ref{prop_fps} that it is given by 
\begin{equation}
\label{mrr9}
\psi_-(s,0) = \frac1\sigma c_-(s,0) \e^{-(2-\deltap)^2/2\sigma^2v_-(s,0)},
\end{equation}
where $\sigma^2v_-(s,0)$ is the variance of $y^-_s$, see~\eqref{mrr2a}, and
the prefactor $c_-(s,0)$ does not play an important r\^ole. Since $v_-(s,0)$
behaves asymptotically like $\vper_-(s)$, the first-passage times at
$1-\deltap$ are sharply concentrated in small neighbourhoods of the local
maxima of $\vper_-(s)$. All these maxima correspond to a single point in
space, but with a different number of revolutions around the stable orbit.
In fact, the theory of large deviations allows to establish a qualitatively
similar behaviour in the general, nonlinear case. 

Next we consider the density $t\mapsto q(t,s)$ of paths starting at time
$s$ in $1-\deltap$ and reaching the unstable orbit in $+1$ at time $t$. In
Sections~\ref{ssec_fpu} and \ref{ssec_fpK}, we establish the expression 
\begin{equation}
\label{mrr10}
q(t,s) = \frac1\sigma \overline c_+(t,s) \e^{-\deltap^2/2\sigma^2\vhat_+(t,s)},
\end{equation}
where $\overline c_+(t,s)$ decays like $\e^{-2\alpha(t,s)}$ and
$\vhat_+(t,s)=\e^{-2\alpha(t,s)}v_+(t,s)$. The exponential decay of the
prefactor is due to the fact that a large, asymptotically constant fraction
of paths leave the neighbourhood of the unstable orbit during each period. 
When computing $q(t,s)$, we have to take into account all paths crossing
the levels $1-\deltap$ and $1-\deltam$ arbitrarily often, before reaching
$+1$. This is done with the help of a renewal equation, discussed in
Section~\ref{sec_re}. The main contribution, however, comes from paths
going directly from $1-\deltap$ to $+1$, without returning to $1-\deltam$. 

For $t\gg s$, $\vhat_+(t,s)$ approaches the function $\vhatper_+(s)$, which
is periodic in the {\em initial\/} time~$s$. This means that paths reaching
$+1$ at some fixed time $t$ have left $1-\deltap$ with approximately equal
probability near any local maximum of $\vhatper_+(s)$. This {\em bottleneck
effect\/} is due to the fact that the level $1-\deltap$ corresponds to
a noncharacteristic boundary where large-deviation results guarantee
concentration of paths, while the unstable orbit causes a strong
dispersion of paths.

The first-passage density $p_+(t)$ is given by the integral of
$q(t,s)\psi_-(s,0)$, which can be evaluated by the Laplace method, yielding a
sum over all minima of the function
\begin{equation}
\label{mrr10.1}
s \mapsto \frac{\deltap^2}{\vhat_+(t,s)} + \frac{(2-\deltap)^2}{v_-(s,0)} .
\end{equation} 
For $0\ll s\ll t$, this function is close to $\rhoper(s)^2$, and thus has
one minimum per period. A careful analysis, given in
Section~\ref{ssec_xll}, shows that for $t$ near $nT$, the $k$th term of the
sum is proportional to 
\begin{equation}
\label{mrr11}
\sigma\biggbrak{
\frac{\gamma(t)}{2\sigma^2}
\exp\Bigset{-2(n-k)\lambda T - \frac{\gamma(t)}{2\sigma^2}
\e^{-2(n-k)\lambda T}}},
\end{equation} 
for some $\gamma(t)$ (the same expression is obtained in \cite{MS4}, using
WKB approximations).  This term is the contribution of the paths making $k$
revolutions around the stable orbit before reaching $1-\deltap$, and $n-k$
revolutions afterwards. Adding $\lambda T$ to $\abs{\log\sigma}$ will
multiply $\sigma^2$ by $\e^{-2\lambda T}$, but this only results in a
rearrangement of the terms in the sum, without changing its value up to
small boundary terms. This roughly explains why the prefactor is periodic
in $\abs{\log\sigma}$. A more detailed analysis of the sum, given in
Section~\ref{ssec_xls}, is necessary to obtain the precise
form~\eqref{mrr6}. 

\begin{figure}
\centerline{\includegraphics*[clip=true,width=125mm]{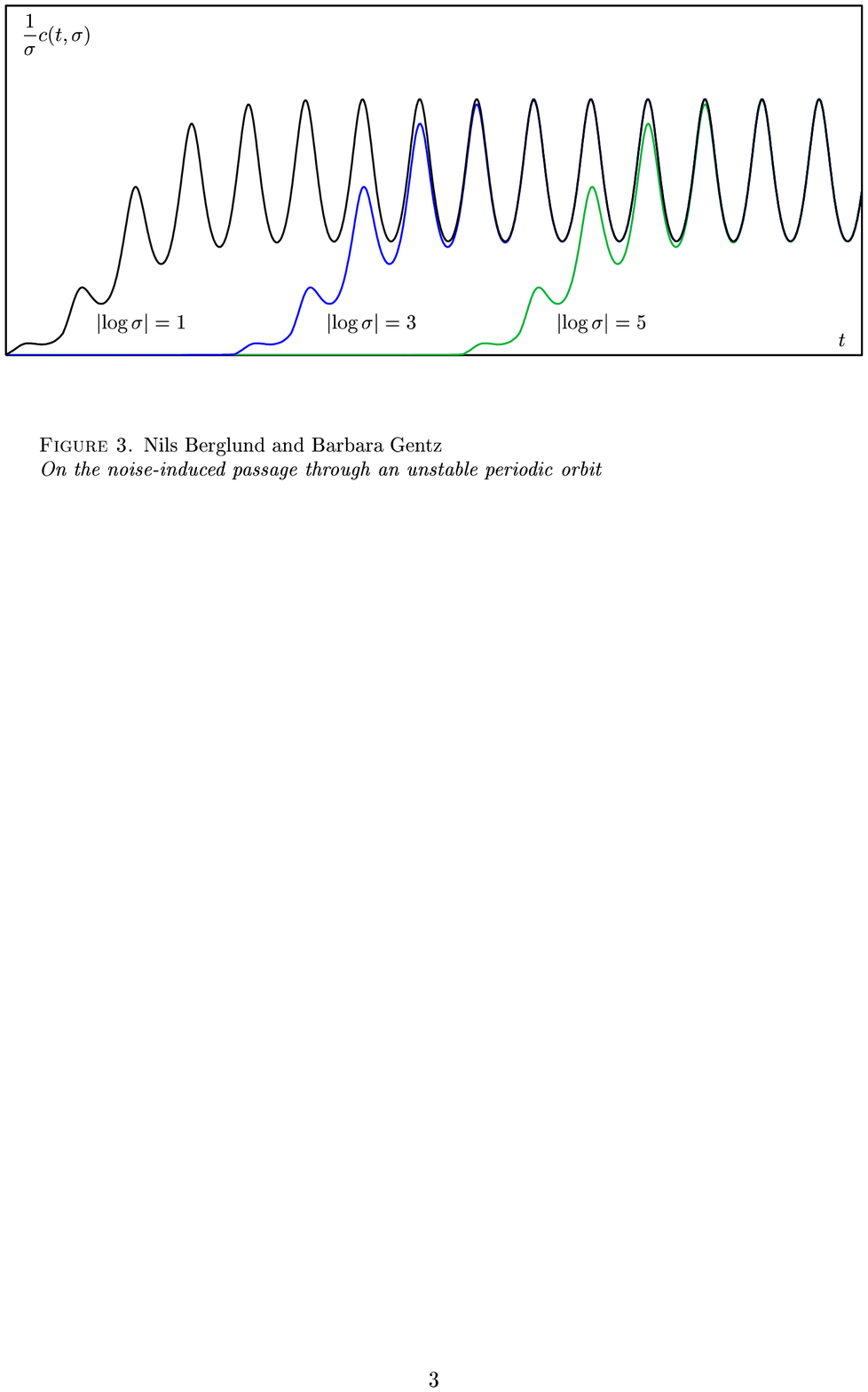}}
 \captionspace
 \caption[]
 {The prefactor $c(t,\sigma)$ of the first-passage density as function of
 time, for different values of $\sigma$. After a  transient regime of
 duration $2\abs{\log\sigma}/\lambda$, the density approaches a periodic
 function of time. Only for times exponentially large in $1/\sigma^2$ will
 the density show a visible decay.}
\label{fig3}
\end{figure}


\subsection{Discussion}
\label{ssec_mrd}

\subsubsection*{Exponential asymptotics}

The exponential rate $R^2$ occurring in~\eqref{mrr4} is independent of $t$,
meaning that on the level of exponential asymptotics, all points on the
unstable orbit are reached with the same probability. This is in sharp
contrast with, say, points on the intermediate level $1-\deltap$, which are
reached with nonconstant exponential rate $(2-\deltap)^2/2v_-(t,0)$. This
difference is a natural consequence of the fact that the unstable orbit is a
characteristic boundary, unlike the level $1-\deltap$. The rate $R^2/2$ 
is the value of the so-called {\em boundary quasipotential\/} of the
Wentzell--Freidlin theory.  

Note that the expressions~\eqref{mrr4} and~\eqref{mrr6} for the
first-passage density are not in contradiction with the property
\begin{equation}
\label{mrd0}
\lim_{\sigma\to0} \bigprob{\e^{(R^2-\delta)/2\sigma^2} \leqs \tau_+ \leqs
\e^{(R^2+\delta)/2\sigma^2}} =1,
\end{equation}
holding according to the Wentzell--Freidlin theory for any fixed $\delta>0$. 
This expression suggests that first-passage times
are concentrated near Kramers' time $\e^{R^2/2\sigma^2}$, but the fact that
$\delta$ is independent of $\sigma$ actually allows for a large
spreading of the distribution of $\tau_+$ which manifests itself in
the behaviour of the prefactor $c(t,s)$. 

\begin{figure}
\centerline{\includegraphics*[clip=true,width=125mm]{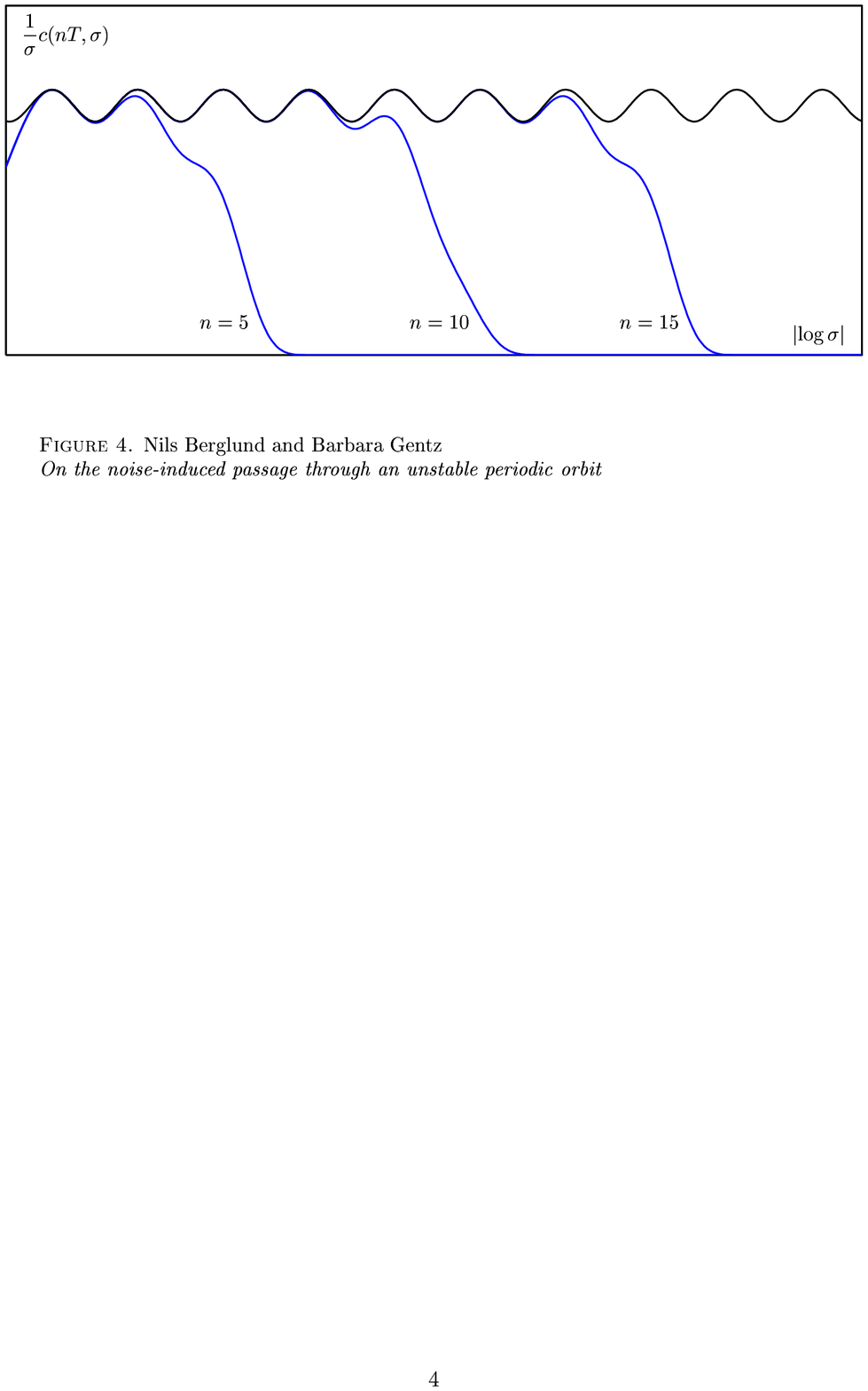}}
 \captionspace
 \caption[]
 {The prefactor of the first-passage density as function of
 $\abs{\log\sigma}$, for three different values of $t=nT$. The prefactor is
 close to a periodic function of $\abs{\log\sigma}$ for
 $2\abs{\log\sigma} \leqs \alpha(nT)=\lambda nT$, and decreases
 exponentially, with an 
 exponent proportional to $-\const \e^{-\alpha(nT)}/\sigma^2$, for
 $2\abs{\log\sigma}> \alpha(nT)= \lambda nT$.}
\label{fig4}
\end{figure}

\subsubsection*{Time scales}

Two time scales play a r\^ole for the behaviour of the first-passage density
$p_+(t)$. They delimit three time domains, whose boundaries, however, are
not particularly sharp. The time scale $\trelax=2\abs{\log\sigma}/\lambda$
is the time needed for the variances $v_-$ and $\vhat_+$ to approach their
asymptotic values; the effect of the transient phase is still visible, for
$t\geqs\trelax$, in the error term $\Order{\e^{-\alpha(t)}/\sigma^2}$. 
The metastable time scale $\tK=\e^{R^2/2\sigma^2}$,
often called Kramers' time, measures roughly the time needed for a
substantial fraction of paths to reach the unstable orbit. There is no sharp
transition, however, between metastable and asymptotic regime, one rather
expects the density $p_+(t)$ to decrease geometrically from one period to
the next one for all times $t\gg\trelax$. 

\subsubsection*{Oscillatory behaviour and cycling}

Consider now the behaviour of the first-passage density $p_+(t)$
during one given 
period $[nT,(n+1)T]$ in the metastable regime. The leading term
in~\eqref{mrr6} depends periodically on $t$ and $\abs{\log\sigma}$. For
fixed noise intensity $\sigma$, $p_+(t)$ is close to a periodic function of
$t$, see~\figref{fig3}. More surprisingly, for fixed $t$ (which also
means a fixed position on the unstable orbit), $p_+$ is proportional
to a term depending periodically on $\abs{\log\sigma}$ (\figref{fig4})
-- this is the phenomenon pointed out in \cite{MS4}. 

\begin{figure}
\centerline{\includegraphics*[clip=true,width=125mm]{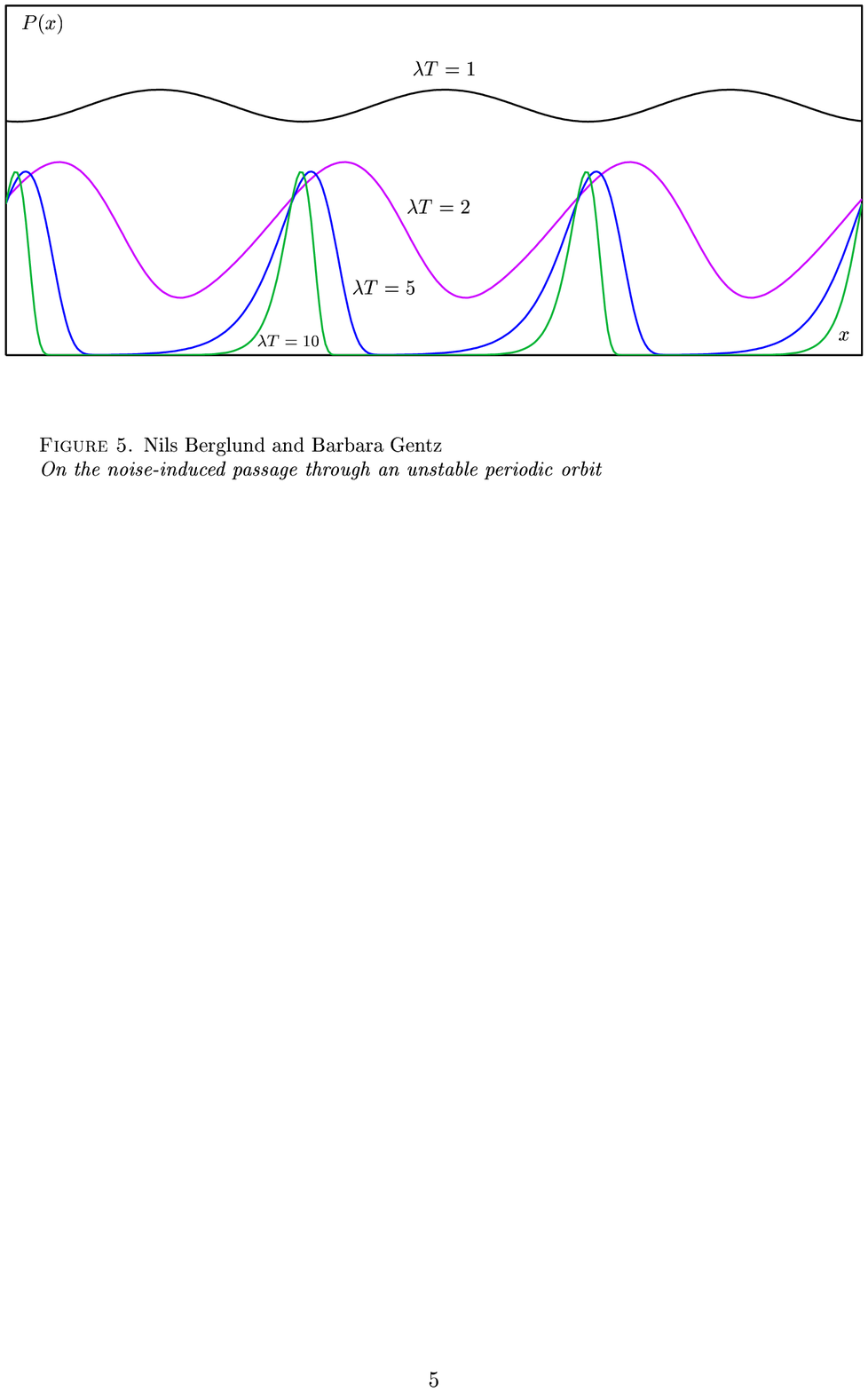}}
 \captionspace
 \caption[]
 {The cycling profile $P(x)$, plotted over three periods, for four different
 values of $\lambda T$. The normalization is such that the integral of
 $P(x)$ over one period is equal to $1/2\lambda T$.}
\label{fig5}
\end{figure}

It is, however, fundamental to consider the joint dependence of $p_+$ on
$t$ and $\sigma$: The leading term in~\eqref{mrr6} can be viewed as a
periodic \lq\lq profile\rq\rq\ function of $\abs{\log\sigma}-\theta(t)$,
modulated by the periodic function $\theta'(t)$. The maximum of $t\mapsto
p_+(t)$ moves once around the unstable orbit as $\abs{\log\sigma}$
increases by $\lambda T$: We thus recover the phenomenon of {\em cycling\/}
described in \cite{Day3,Day4} for the first-exit location, but in addition
we obtain here the same behaviour for the first-passage time, which keeps
track of the winding number around the unstable orbit. Also note that the
probability of reaching the unstable orbit during $[nT,(n+1)T]$ is
proportional to 
\begin{equation}
\label{mrd1}
\int_{nT}^{(n+1)T} \theta'(t)
P\Bigpar{\frac{\abs{\log\sigma}-\theta(t)}{\lambda T}} \6t 
= \lambda T\int_{x_0-1}^{x_0} P(x) \6x = \frac12,
\end{equation}
where we have set $x_0=(\abs{\log\sigma}-\theta(nT))/\lambda T$. While
the peak of $p_+(t)$ moves around the unstable orbit as $\sigma$ decreases,
changing its height periodically, the area below $p_+(t)$ remains constant.
Another way to interpret~\eqref{mrr6} is to consider $\theta(t)$ as 
intrinsic time: the first-passage density expressed with respect to the
time $\theta(t)$ is translated around the unstable orbit as $\sigma$
decreases, with constant \lq\lq velocity\rq\rq, without changing its
shape. 

\subsubsection*{Bottleneck effect}

Paths reach the intermediate level $1-\delta_1$ at times 
concentrated in small windows around $s^\star + kT$. When approaching the
unstable orbit, they are strongly dispersed. As a result, a path reaching
$y=+1$ at a fixed time $t\in[nT,(n+1)T)$ may have idled along the unstable
orbit for an almost arbitrarily long time span. The probability that it has
come through the $k$th window is proportional to~\eqref{mrr11}. When
translating this from time to space, all windows become superimposed in
$\ph^\star=(2\pi/T)s^\star$, meaning that with high probability, all paths
reaching the unstable orbit have crossed the curve $y=1-\delta_1$ through
the same small bottleneck. In fact, the same is true for the crossing of
any curve bounded away from the periodic orbits. We thus recover the
well-known fact that the transition between the orbits is likely to occur
in a small neighbourhood of a fixed trajectory spiraling away from the
stable orbit, the so-called most probable exit path. But note that
the concept of most probable exit paths becomes irrelevant for the
dynamics close to the unstable orbit.

\subsubsection*{High-frequency limit $T\ll 1$}

We now examine how the first-passage density changes as a function of the
period $T$. By this we mean that the functions $a$ and $g$ are scaled by a
factor $T$, that is, $a(t)=a_T(t)=a_1(t/T)$ and $g(t)=g_T(t)=g_1(t/T)$ for
some fixed functions $a_1$ and $g_1$. In particular,
$\alpha(t)=T\alpha_1(t/T)$, so that, for instance, \eqref{mrr1a} becomes 
\begin{equation}
\label{mrd2}
\vper_-(Ts) = \frac{T}{1-\e^{-2\lambda T}}
\int_s^{s+1} \e^{-2T\alpha_1(s+1,u)}g_1(u)^2\6u. 
\end{equation}
In the limit $T\to0$, the functions $\vper_-$ and $\vhatper_+$ both
approach a constant value $\average{g_1^2}/2\lambda$, where
$\average{g_1^2}$ denotes the average value of $g_1^2$. The exponential
rate $R^2$ approaches the value
$2\lambda\brak{\deltap^2+(2-\deltap)^2}/\average{g_1^2}$ (of course, we
cannot actually take the limit $T\to0$, because $\rhoper$ would also become
constant in this limit, so that the discussion is to be understood for small
but finite $T$). 

Since $\abs{\Gamma(1-\icx y)}$ decreases exponentially with $\abs{y}$, 
$P(x)$ is close to a sinuso\"{\i}d with amplitude of order
$\e^{-\const/T}$, so that the cycling profile becomes flat in the
high-frequency limit, cf.~\figref{fig5}. This is related to the fact that many
terms, i.\,e., paths with many different winding numbers, contribute to the
sum~\eqref{mrr8}. In the limit, the cycling velocity $1/\th'(t)$ behaves like 
$\average{g_1^2}/\lambda g(t)^2$, so that the first-passage density
$p_+(t)$ behaves like the noise coefficient $g(t)^2$, independently of
$a(t)$. 

\subsubsection*{Adiabatic limit $T\gg 1$}

In the low-frequency limit, $\vper_-(t)$ and $\vhatper_+(t)$ both follow
adiabatically $v^\star(t)$, at a distance of order $1/T$. As a result, the
exponential rate $R^2$ approaches $\brak{\deltap^2+(2-\deltap)^2}/\bv$, which
is smaller than its value in the limit $T\to0$ (because $2\lambda\bv >
\average{2a_1v^\star}=\average{g_1^2}$). This is due to the fact that paths
have enough time to probe for the moment when the transition costs the
least. 

The sum~\eqref{mrr8} is dominated by one term, so that the cycling profile
$P(x)$ is sharply peaked, cf.~\figref{fig5}. The first-passage density is thus
dominated by paths making a fixed number of revolutions around the unstable
orbit.  The cycling velocity $1/\th'(t)$ converges to $1/a(t)$. The
first-passage density $p_+(t)$ thus consists of a peak of height
proportional to $a(t)$, moving around with velocity $1/a(t)$. 

In spite of the fact that $P(x)$ is sharply peaked, one can show that the
first-passage density at fixed $t$ is still a monotonously decreasing
function of $\sigma$. This is due to the exponentially small factor
$\e^{-R^2/2\sigma^2}$. 



\goodbreak
\section{The renewal equation}
\label{sec_re}


In this section, we establish a renewal equation satisfied by the
first-passage time $\tau_+$ at $+1$ of the process $y_t$  defined in
Section~\ref{ssec_mrs}. By restarting the process at the first time $\tauup
= \tauup(0,-1)$ at which the level $1-\deltap$ is reached, the distribution
function of $\tau_+$ can be written as 
\begin{align}
\nonumber
\probin{0,-1}{\tau_+ \leqs t} 
&= \bigexpecin{0,-1}{\indexfct{\tauup \leqs t}
  \probin{\tauup,1-\deltap}{\tau_+ \leqs t}}
\\
\label{re1}
&= \int_0^t Q(t,s) \psi_-(s,0) \6s,
\end{align}
where we have introduced the quantities 
\begin{align}
\label{re2a}
Q(t,s) &= \probin{s,1-\deltap}{\tau_+ \leqs t} \\
\label{re2b}
\psi_-(s,0) &= \dpar{}{s}\probin{0,-1}{\tauup \leqs s}.  
\end{align}
The first-passage density $\psi_-(s,0)$ of $y_s$ to $1-\deltap$ depends only
on the linear process $y^-_s$; we will discuss its computation in
Section~\ref{ssec_fps}. The function $Q(t,s)$ depends on all subsequent
switchings of the process $y_t$ between $y_t^+$ and $y_t^-$, and therefore
we will write it as the solution of an integral equation, that will
serve as a {\it renewal equation}. It is obtained by restarting the
process each time the level $1-\deltap$ is reached from below.  

Let $y^\#_t$ denote the stochastic process obtained by killing $y^+_t$
upon first hitting the level $+1$, and introduce the stopping time  
\begin{equation}
\label{re3}
\taudown^{\#}(t_0) = \inf\setsuch{s>t_0}{y^\#_s < 1-\deltam}, 
\end{equation}
with the convention that $\taudown^{\#}=\infty$ if $y^+_s$ either hits
level~$+1$ before reaching $1-\deltam$ or never reaches $1-\deltam$.  

\begin{prop}
\label{prop_re}
$Q(t,s)$ satisfies the\/ {\em renewal equation}
\begin{equation}
\label{re4}
Q(t,s) = P_1(t,s) + \int_s^t Q(t,u) K(u,s) \6u,
\end{equation}
where 
\begin{align}
\label{re5a}
P_1(t,s) &= \probin{s,1-\deltap}{\tau_+ \leqs t,\tau_+<\taudown}, \\
\label{re5b}
K(u,s) &= \int_s^u \psiup(u,v) \psidown^{\#}(v,s) \6v, \\
\label{re5c}
\psiup(u,v) &= \dpar{}{u} \probin{v,1-\deltam}{\tauup \leqs u}, \\
\label{re5d}
\psidown^{\#}(v,s) &= \dpar{}{v} \probin{s,1-\deltap}{\taudown^{\#}
  \leqs v}. 
\end{align}
\end{prop}
\begin{proof}
Splitting up the event according to whether $\taudown>\tau_+$ or
$\taudown\leqs\tau_+$, we can write 
\begin{equation}
\label{re6:1}
Q(t,s) = P_1(t,s) + \probin{s,1-\deltap}{\taudown\leqs\tau_+ \leqs t}.  
\end{equation}
Note that $\taudown\leqs\tau_+ \leqs t$ if and only if the killed process
$y^{\#}_t$ first reaches $1-\deltam$ and then $+1$, and both events
occur before time $t$. Therefore, the second term 
on the right-hand side of~\eqref{re6:1} can be written as 
\begin{align}
\nonumber
&\Bigexpecin{s,1-\deltap}{\indexfct{\taudown^{\#} \leqs t} 
\bigprobin{\taudown^{\#},1-\deltam}{\tau_+ \leqs t}}\\
\nonumber
&\qquad = \Bigexpecin{s,1-\deltap}{\indexfct{\taudown^{\#} \leqs t} 
\bigexpecin{\taudown^{\#},1-\deltam}
{\indexfct{\tauup \leqs t}
Q(t,\tauup)}}\\
&\qquad = \int_s^t \int_v^t Q(t,u) \psiup(u,v) \6u \mskip6mu
\psidown^{\#}(v,s) \6v. 
\label{re6:2}
\end{align}
Using Fubini's theorem to interchange the integrals, we obtain the second
term in~\eqref{re4}. 
\end{proof}

Note that since $Q(t,t)=0$, the derivatives 
\begin{equation}
\label{re7}
q(t,s) = \dpar{}{t}Q(t,s), 
\qquad\qquad
p_1(t,s) = \dpar{}{t}P_1(t,s)
\end{equation}
satisfy the similar renewal equation 
\begin{equation}
\label{re8}
q(t,s) = p_1(t,s) + \int_s^t q(t,u) K(u,s) \6u. 
\end{equation}

The renewal equation~\eqref{re4} can be solved by iterations. In our case,
$K$ will be so small that only the first term needs to be computed. The
convergence of the iterative method can, however, easily be controlled. 

\begin{prop}
\label{prop_re2}
Let $K_1(u,s)=K(u,s)$ and define for $n\in\N$ 
\begin{equation}
\label{re9}
K_{n+1}(u,s) = \int_s^u 
K(u,v) K_n(v,s) 
\6v. 
\end{equation}
Then, for any $N\geqs1$, the solution of the renewal equation~\eqref{re4}
satisfies 
\begin{equation}
\label{re10}
Q(t,s) = P_1(t,s) + \sum_{n=1}^N \int_s^t P_1(t,u)K_n(u,s)\6u +
R_{N+1}(t,s), 
\end{equation}
with a remainder satisfying $0\leqs R_N(t,s)\leqs
M^N\abs{t-s}^N/N!$ for some constant~$M$ which depends
neither on $N$ nor on $t,s$. 
Hence, $Q(t,s)$ can be written as a series which converges uniformly
on compact sets.   
\end{prop}
\begin{proof}
First note that $\psiup(t,u)$ is bounded by some constant~$M>0$ which 
implies that $K(t,s)\leqs M$ for all $t,s$. By induction, we see
that~\eqref{re10} holds with  
\begin{equation}
\label{re11a}
R_N(t,s) = \int_s^t Q(t,v)K_N(v,s) \6v
\end{equation}
and that 
\begin{equation}
\label{re11b}
\int_s^t K_N(v,s) \6v \leqs M^N\abs{t-s}^N/N!
\end{equation}
Thus the bounds on $R_N(t,s)$ follow from the trivial bounds $0\leqs Q(t,v)
\leqs 1$. 
\end{proof}

The $n$th term in the sum in~\eqref{re10} is the contribution of those
paths which cross the levels $1-\deltap$ and $1-\deltam$ alternatively
$n$ times during the time interval $[s,t]$, before finally reaching~$+1$. 



\section{The first-passage densities}
\label{sec_fp}


\subsection{Leaving the stable orbit}
\label{ssec_fps}

\begin{figure}
\centerline{\includegraphics*[clip=true,width=100mm]{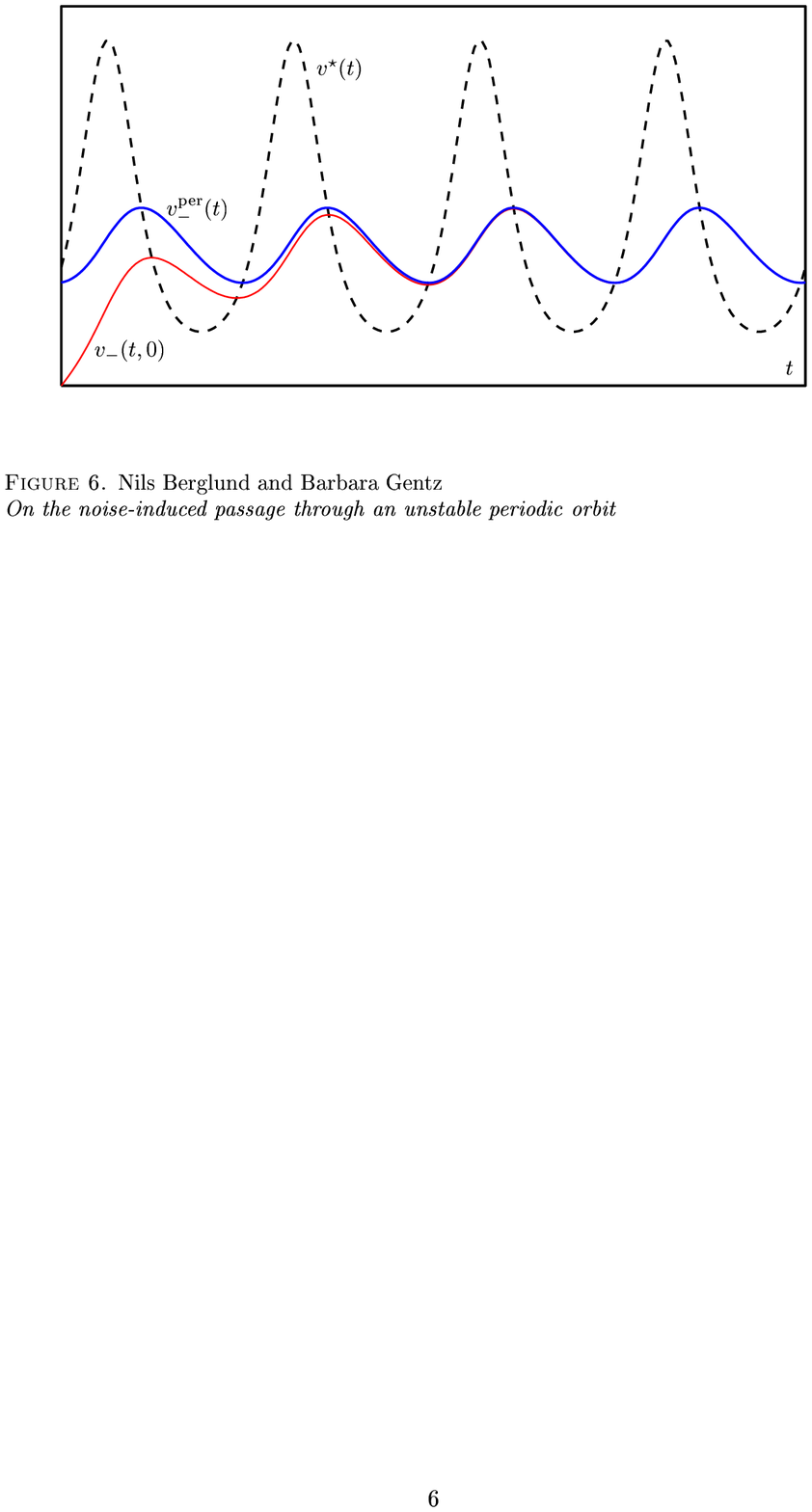}}
 \captionspace
 \caption[]
 {The functions $v_-(t,0)$ and $\vper_-(t)$.}
\label{fig6}
\end{figure}

In this section, we ignore all possible switchings between
the processes $y^+_t$ and $y^-_t$ in the definition of $y_t$ and focus
on the stochastic process $y^-_t(t_0,-1)$ only. Recall that this
stochastic process is defined by the SDE 
\begin{equation}
\label{fps1}
\6 y^-_t = - a(t) (y^-_t+1) \6t + \sigma g(t) \6W_t,
\qquad
y^-_{t_0}=-1,
\end{equation}
for some initial time $t_0$. Thus it is a Gaussian process, given by 
\begin{equation}
\label{fps2}
y_t = y^-_t(t_0,-1) = -1 + \sigma\int_{t_0}^t \e^{-\alpha(t,s)}g(s)\6W_s,
\end{equation}
where $\alpha(t,s) = \alpha(t)-\alpha(s) = \int_s^t a(u)\6u$. At
time~$t$, the Gaussian random variable $y^-_t$ has mean~$-1$ and
variance $\sigma^2 v_-(t,t_0)$, where  
\begin{equation}
\label{fps3}
v_-(t,t_0) = \int_{t_0}^t \e^{-2\alpha(t,s)} g(s)^2 \6s. 
\end{equation}
Note that $v_-(t,t_0)$ satisfies the deterministic differential equation 
\begin{equation}
\label{fps4}
\dpar{}{t} v_-(t,t_0) = -2a(t) v_-(t,t_0) + g(t)^2, 
\end{equation}
with initial condition $v_-(t_0,t_0)=0$. This equation also
admits a periodic solution
\begin{equation}
\label{fps5}
\vper_-(t) = \frac1{1-\e^{-2\lambda T}} \int_t^{t+T} \e^{-2\alpha(t+T,s)}
g(s)^2\6s.
\end{equation}
Recall that $\lambda = \alpha(T)/T$ is the Lyapunov exponent. These two
solutions of~\eqref{fps4} are related by 
\begin{equation}
\label{fps6}
v_-(t,t_0) = \vper_-(t) - \e^{-2\alpha(t,t_0)}\vper_-(t_0), 
\end{equation}
see~\figref{fig6}. 

The following result describes the behaviour of the first-passage density
$\psi_-(t,0)$ of $y^-_t(0,-1)$ to $1-\deltap$.

\begin{prop}
\label{prop_fps}
Assume that $\vper_-(t)\leqs 2v^\star(t)(1-\Delta)$ for some constant
$\Delta>0$. Then the first-passage density $\psi_-(t,0)$ can be written as 
\begin{equation}
\label{fps8}
\psi_-(t,0) = \frac1\sigma c_-(t,0) \e^{-\rho_-(t,0)^2/2\sigma^2}, 
\end{equation}
where 
\begin{equation}
\label{fps9}
\rho_-(t,0)^2 = \frac{(2-\deltap)^2}{v_-(t,0)}
\end{equation}
and
\begin{equation}
\label{fps10}
c_-(t,0) = \frac{2-\deltap}{\sqrt{2\pi}}
\biggbrak{\frac{1}{v_-(t,0)} - \frac1{2v^\star(t)}}
\frac{g(t)^2}{\sqrt{v_-(t,0)}}
\biggbrak{1+\frac1\Delta\biggOrder{\frac\sigma{\Delta^2} + 
\frac{\e^{-\const\Delta^2/\sigma^2}}{\sigma} t}}
\end{equation}
for all $t < \tfrac12\sigma(1-\sigma/\Delta^2) \e^{\const\Delta^2/\sigma^2}$.
\end{prop}
\begin{proof}
The process $z^-_t = \e^{\alpha(t)}(y^-_t+1)$ satisfies a stochastic
differential equation without drift term, 
\begin{equation}
\label{fps11}
\6z^-_t = \sigma\e^{\alpha(t)}g(t)\6W_t,
\qquad z^-_0=0.
\end{equation}
It is a Gaussian process, with mean zero and variance
$\sigma^2\vtilde_-(t,0)=\sigma^2\e^{2\alpha(t)}v_-(t,0)$. The first passage
of $y^-_t$ at $1-\deltap$ corresponds to the first passage of $z^-_t$ at the
time-dependent level $d_-(t) = (2-\deltap)\e^{\alpha(t)}$. Some properties
of first-passage densities of such Gaussian processes are discussed in
the appendix.

We will apply Corollary~\ref{cor_appendix} from the appendix, after checking
that $\vtilde_-(t,0)$ and $d_-(t)$ satisfy the required
conditions. \eqref{app5a} can easily be verified by rewriting the condition
with the help of $\vper_-(t)$ and recalling the assumption on
the relation between $\vper_-(t)$ and $v^\star(t)$. Note that $\Delta$
in this assumption and $\Delta$ in~\eqref{app5a} differ by a
constant. Checking~\eqref{app5c} is even easier, just keep in mind
that $v_-(t,0)$ approaches $\vper_-(t)$ and thus is bounded. 
Next,~\eqref{app5b} can be established by showing first that there
exists a $\delta_0>0$ such that $\abs{\tilde c(t,s)}$ as defined in the
appendix is bounded above by $\const \sqrt{t-s}$ as long as 
$t-s \leqs \delta_0$. For $t-s > \delta_0$, one can show that
$\abs{\tilde c(t,s)}$ is at most of order one. Finally,
Assumption~\eqref{app9.1} is seen to be satisfied for large enough
$M_3$ by comparing again $v_-(t,0)$ with $\vper_-(t)$. Applying the
corollary immediately yields~\eqref{fps8}--\eqref{fps10}.
\end{proof}


\subsection{Reaching the unstable orbit}
\label{ssec_fpu}

We now turn to the process $y^+_t(s,1-\deltap)$, defined by the SDE 
\begin{equation}
\label{fpu1}
\6 y^+_t = a(t) (y^+_t-1) \6t + \sigma g(t) \6W_t,
\qquad
 y^+_{s}=1-\deltap,
\end{equation}
for some initial time $s$. It is given by 
\begin{equation}
\label{fpu2}
y^+_t = y^+_t(s,1-\deltap) 
= 1 - \deltap\e^{\alpha(t,s)} + \sigma\int_{s}^t \e^{\alpha(t,u)}g(u)\6W_u.
\end{equation}
At time~$t$, $y^+_t$ is Gaussian with variance $\sigma^2 v_+(t,s)$,
where 
\begin{equation}
\label{fpu3}
v_+(t,s) = \int_s^t \e^{2\alpha(t,u)} g(u)^2 \6u. 
\end{equation}
As $y^+_t$ is spreading fast as $t$ increases, it is helpful to
consider the stochastic process
\begin{equation}
\label{fpu3.1}
z^+_t = \deltap + \e^{-\alpha(t,s)}[y^+_t-1],
\end{equation}
which is the solution of the SDE
\begin{equation}
\label{fpu3.2}
\6 z^+_t = \sigma \e^{-\alpha(t,s)}g(t) \6W_t,
\qquad
z^+_s = 0.
\end{equation}
It is also Gaussian and has variance $\sigma^2 \vhat_+(t,s)$, where
\begin{equation}
\label{fpu4}
\vhat_+(t,s) = \e^{-2\alpha(t,s)}v_+(t,s)
= \int_s^t \e^{-2\alpha(u,s)} g(u)^2 \6u.
\end{equation}
Note that $s\mapsto \vhat_+(t,s)$ satisfies the deterministic
differential equation  
\begin{equation}
\label{fpu5}
\dpar{}{s} \vhat_+(t,s) = 2a(s) \vhat_+(t,s) - g(s)^2,
\end{equation}
with $\vhat_+(t,t)=0$. This equation also admits a
periodic solution 
\begin{equation}
\label{fpu6}
\vhatper_+(s) = \frac1{\e^{2\lambda T}-1} \int_s^{s+T} \e^{2\alpha(s+T,u)}
g(u)^2\6u.
\end{equation}
The function $\vhat_+(t,s)$ can then be determined by the relation
\begin{equation}
\label{fpu7}
\vhat_+(t,s) = \vhatper_+(s) - \e^{-2\alpha(t,s)}\vhatper_+(t),
\end{equation}
see~\figref{fig6}.
By the reflection principle, the distribution function of the first-passage
time~$\tildetau_+$ of $y^+_t$ at $+1$ is given by 
\begin{equation}
\label{fpu8}
\probin{s,1-\deltap}{\tildetau_+ \leqs t} 
= 2 \probin{s,1-\deltap}{y^+_t \geqs 1}
= 2 \Phi\Bigpar{-\frac{\rho_+(t,s)}\sigma},
\end{equation}
where $\Phi(x) = (2\pi)^{-1/2}\int_{-\infty}^x \e^{-u^2/2}\6u$ denotes the
distribution function of the standard normal law, and 
\begin{equation}
\label{fpu9}
\rho_+(t,s)^2 = \frac{\deltap^2 \e^{2\alpha(t,s)}}{v_+(t,s)} 
= \frac{\deltap^2}{\vhat_+(t,s)}. 
\end{equation}
The density of $\tildetau_+$ can thus be written as 
\begin{align}
\nonumber
\dpar{}{t}\probin{s,1-\deltap}{\tildetau_+ \leqs t} 
&= \frac1\sigma \frac{\deltap}{\sqrt{2\pi}}
\frac{\partial_t \vhat_+(t,s)}{\vhat_+(t,s)^{3/2}} 
\e^{-\rho_+(t,s)^2/2\sigma^2} \\
&=\frac1\sigma \frac{\deltap}{\sqrt{2\pi}}
\frac{g(t)^2 \e^{-2\alpha(t,s)}}
{\vhat_+(t,s)^{3/2}} 
\e^{-\rho_+(t,s)^2/2\sigma^2}. 
\label{fpu10}
\end{align}
Let us now show that the function 
\begin{equation}
\label{fpu11}
p_1(t,s) 
= \dpar{}{t} \probin{s,1-\deltap}{\tau_+ \leqs t,\tau_+<\taudown}
= \dpar{}{t} \probin{s,1-\deltap}{\tildetau_+ \leqs t,\tildetau_+<\taudown}
\end{equation}
behaves in the same way. Recall that $\taudown=\taudown(s,1-\deltap)$
denotes the first-passage time of $y^+_t(s,1-\deltap)$ at the lower
level $1-\deltam$. It is obvious that $p_1(t,s)$ is bounded above
by~\eqref{fpu10},  but we want to show that it actually has the same
exponential asymptotics, and, moreover, almost the same prefactor. In other
words, it is not only unlikely that a path reaches the unstable orbit, but
also the conditional probability that a path returns to  $1-\deltam$ before
reaching~$+1$, given that it actually reaches the unstable orbit at~$+1$,
is small. Hypothesis H3 plays a crucial r\^ole here. 

In a first step, we study the density of $\taudown$.

\begin{lemma}
\label{lem_fpu}
The density $\psidown(u,s) = \dpar{}{u}\probin{s,1-\deltap}{\taudown
  \leqs u}$ of $\taudown = \taudown(s,1-\deltap)$ satisfies 
\begin{equation}
\label{fpu12}
\psidown(u,s) = \frac1\sigma \cdown(u,s)
\e^{-\rhodown(u,s)^2/2\sigma^2}, 
\end{equation}
where the prefactor $\cdown(u,s)$ is bounded by a constant times
$\vhat_+(u,s)^{-3/2}\e^{-\alpha(u,s)}$ and
\begin{equation}
\label{fpu13}
\rhodown(u,s)^2 = 
\frac{\bigbrak{\deltap-\deltam\e^{-\alpha(u,s)}}^2}{\vhat_+(u,s)}.
\end{equation}
\end{lemma}
\begin{proof}
Recall the definition of the stochastic process $z^+_u = \deltap +
\e^{-\alpha(u,s)}[y^+_u-1]$ from~\eqref{fpu3.1}. The first passage of
$y^+_u$ at $1-\deltam$ coincides with the first passage of $z^+_u$ at the
time-dependent level $d_+(u) = \deltap - \deltam
\e^{-\alpha(u,s)}$. Therefore, the result follows from
Lemma~\ref{lem_appendix0} of the appendix which we apply with
$d(u)=d_+(u)$ and $v(u) = \vhat_+(u,s)$ for fixed~$s$. 
\end{proof}

\begin{prop}
\label{prop_fpu}
Assume that $\vhatper_+(u) \leqs 2 v^\star(u)(1-\Delta)$ for all~$u$. Then  
\begin{equation}
\label{fpu16}
p_1(t,s) = \frac1\sigma c_+(t,s) \e^{-\rho_+(t,s)^2/2\sigma^2}, 
\end{equation}
where $\rho_+(t,s)$ is the function defined in~\eqref{fpu9} and 
\begin{equation}
\label{fpu17}
c_+(t,s) = \frac{\deltap}{\sqrt{2\pi}}
\frac{g(t)^2 \e^{-2\alpha(t,s)}}{\vhat_+(t,s)^{3/2}} 
 - \frac1\sigma \bigOrder{\e^{-\alpha(t,s)} 
\e^{-\const\Delta_0^2/\sigma^2}},
\end{equation}
where $\Delta_0 = \Delta/(1-\Delta)\wedge(\deltam-\deltap)/\deltap
\wedge 1$.
If the condition on $\vhatper_+$ is not satisfied, \eqref{fpu16} remains
valid, but the prefactor can be smaller than~\eqref{fpu17}. 
\end{prop}
\begin{proof}
The main idea of the proof is that the condition on $\vhatper_+$ excludes
that the most probable path going from $(s,1-\delta_1)$ to $(t,1)$ crosses
the lower level $1-\delta_2$. The only difficulty resides in exploiting this
fact in order to obtain the exponentially small error bound
in~\eqref{fpu17}. 

By definition of $P_1(t,s)$, 
\begin{equation}
\label{fpu18:1}
P_1(t,s) = \probin{s,1-\deltap}{\tildetau_+ \leqs t} 
- \probin{s,1-\deltap}{\taudown<\tildetau_+ \leqs t}.
\end{equation}
Recall that $\taudown^\#$ denotes the first-passage time at level
$1-\deltam$ for the process $y^\#_t$, which is obtained by 
killing $y^+_t$ upon reaching level $+1$. The second term on the
right-hand side can be written as  
\begin{align}
\nonumber
\probin{s,1-\deltap}{\taudown<\tildetau_+ \leqs t} &=
\Bigexpecin{s,1-\deltap}{\indexfct{\taudown^{\#} \leqs t} 
\probin{\taudown^{\#},1-\deltam}{\tildetau_+ \leqs t}}\\
&= \int_s^t \probin{u,1-\deltam}{\tildetau_+ \leqs t} \psidown^\#(u,s) \6u. 
\label{fpu18:2}
\end{align}
Since $\probin{t,1-\deltam}{\tildetau_+ \leqs t}=0$ and
$\psidown^\#(u,s)\leqs\psidown(u,s)$, we have 
\begin{equation}
\label{fpu18:3}
\dpar{}{t} \probin{s,1-\deltap}{\taudown<\tildetau_+ \leqs t} \leqs 
\int_s^t \dpar{}{t}\probin{u,1-\deltam}{\tildetau_+ \leqs t}\psidown(u,s)\6u. 
\end{equation}
The first term in the integrand is similar to~\eqref{fpu10}, the only
difference lying in the initial condition. We can write it in the form 
\begin{equation}
\label{fpu18:4}
\dpar{}{t}\probin{u,1-\deltam}{\tildetau_+ \leqs t} = \frac1\sigma \tilde c_+(t,u)
\e^{-\tilde\rho_+(t,u)^2/2\sigma^2},
\end{equation}
where $\tilde\rho_+(t,u)^2 = \deltam^2/\vhat_+(t,u)$ and $\tilde
c_+(t,u) = \deltam g(t)^2
\e^{-2\alpha(t,u)}/{\sqrt{2\pi}}{\vhat_+(t,u)^{3/2}}$. The derivative  
$p_1(t,s)$ thus satisfies 
\begin{equation}
\label{fpu18:4b}
p_1(t,s) \geqs \frac1\sigma \e^{-\rho_+(t,s)^2/2\sigma^2} 
\biggbrak{\frac{\deltap}{\sqrt{2\pi}}
\frac{g(t)^2 \e^{-2\alpha(t,s)}}{\vhat_+(t,s)^{3/2}}  - \frac1\sigma\int_s^t \tilde c_+(t,u)\cdown(u,s)
\e^{-\chi(u)\rho_+(t,s)^2/2\sigma^2}\6u},
\end{equation}
where 
\begin{equation}
\label{fpu18:5}
\chi(u) = \frac{\tilde\rho_+(t,u)^2 + \rhodown(u,s)^2}{\rho_+(t,s)^2} - 1.
\end{equation}
Assume that $t-s>2$. (For $t-s\leqs2$, it suffices to apply the first
one of the arguments below.) We split the integral
in~\eqref{fpu18:4b} at $u_1 = s+1 $ and 
$u_2 = t - 1$.  The integrals over $[s,u_1]$ and $[u_2,t]$ are
easily seen to be exponentially small because $\chi(u)$ diverges as $u\to
s$ and $u\to t$ (in fact, for small enough $\sigma$, the integrand is
maximal at $u_1$ or $u_2$, respectively). In order to bound the remaining
integral over $[u_1,u_2]$, we only need to find a positive lower bound for
the function $\chi(u)$,  valid whenever $s\leqs u\leqs t$. Using the fact
that  $\vhat_+(t,u) = \e^{2\alpha(u,s)}\brak{\vhat_+(t,s)-\vhat_+(u,s)}$,
it is straightforward to show that 
\begin{equation}
\label{fpu18:6}
\chi(u) = \frac{\brak{\kappa(u)-(1-r(u))}^2}{r(u)(1-r(u))},
\end{equation}
where 
\begin{equation}
\label{fpu18:7}
\kappa(u) = \frac{\deltam}{\deltap}\e^{-\alpha(u,s)}
\qquad\text{and}\qquad
r(u) = \frac{\vhat_+(u,s)}{\vhat_+(t,s)} \in [0,1]. 
\end{equation}
We may assume the existence of a constant $c^\star>1$ such that 
$\vhat_+(t,u)/\vhat_+(t,s)\leqs c^\star$ for all $u\in[s,t]$. 
The function $r(u)$ being increasing with range $[0,1]$, we can define
$u^\star\in(s,t)$ by 
\begin{equation}
\label{fpu18:8}
1 - r(u^\star) = \e^{-2\alpha(u^\star,s)}
\frac{\vhat_+(t,u^\star)}{\vhat_+(t,s)} = \frac1{c^\star}.
\end{equation}
Consider first the case $u\geqs u^\star$. Then 
\begin{equation}
\label{fpu18:9}
\frac{\kappa(u)^2}{1-r(u)} = \biggpar{\frac{\deltam}{\deltap}}^2 
\frac{\vhat_+(t,s)}{\vhat_+(t,u)} 
\geqs \biggpar{\frac{\deltam}{\deltap}}^2\frac1{c^\star}
\geqs \biggpar{\frac{\deltam}{\deltap}}^2(1-r(u)),
\end{equation}
and thus 
\begin{equation}
\label{fpu18:10}
\chi(u) \geqs\frac{\brak{\kappa(u)-(1-r(u))}^2}{1-r(u)}
\geqs \frac1{c^\star} \biggpar{\frac{\deltam}{\deltap}-1}^2
= \frac{(\deltam -\deltap)^2}{c^\star\deltap^2}.
\end{equation}
We turn now to the case $u\leqs u^\star$, in which it suffices to find a lower
bound for $\theta(u)=\kappa(u)+r(u)-1$, since $\chi(u)\geqs 4\theta(u)^2$ for
all $u$. First note that $\theta(s)=(\deltam -
\deltap)/\deltap$ and, since~\eqref{fpu18:8} implies 
$\e^{-\alpha(u^\star,s)}\geqs1/c^\star$, we also know that 
$\theta(u^\star)\geqs(\deltam - \deltap)/c^\star\deltap$. Thus
if $\theta(u)$ reaches its minimum on the boundary of $[s,u^\star]$,
we are done. Otherwise, the fact that 
\begin{equation}
\label{fpu18:11}
\theta'(u) = 0 
\qquad\Leftrightarrow\qquad
\kappa(u) = 2v^\star(u)\frac{\e^{-2\alpha(u,s)}}{\vhat_+(t,s)}
\end{equation}
shows that if $\theta(u)$ reaches its minimum in $(s,t)$, then 
\begin{equation}
\label{fpu18:12}
\theta(u) =
\frac{\e^{-2\alpha(u,s)}}{\vhat_+(t,s)}\bigbrak{2v^\star(u)-\vhat_+(t,u)} 
\geqs \frac1{c^\star}\biggbrak{\frac{2v^\star(u)}{\vhat_+(t,u)}-1}
\geqs \frac{\Delta}{c^\star(1-\Delta)},
\end{equation}
where the last estimate holds due to our assumption on
$\vhatper_+(u)$. Note that if $\vhatper_+(u)$ is not bounded away
from $2v^\star(u)$, then the minimum of $\theta$ is still positive,
but may be of the order $\e^{-2\alpha(t,s)}$, which can become very small. 
\end{proof}


\subsection{The renewal kernel}
\label{ssec_fpK}

The kernel $K(u,s)$ involves the functions $\psiup(u,v)$ and
$\psidown^\#(v,s)$ defined in~\eqref{re5c} and~\eqref{re5d}. We already know
from Lemma~\ref{lem_fpu} that the rate associated with
$\psidown(v,s)$, which provides an upper bound for $\psidown^\#(v,s)$,
is given by  
\begin{equation}
\label{fpK1}
\rhodown(v,s)^2 = 
\frac{\bigbrak{\deltap-\deltam\e^{-\alpha(v,s)}}^2}{\vhat_+(v,s)}.
\end{equation}
Note that $\rhodown(v,s)$ vanishes at a time $v \geqs s$ such that 
$\e^{-\alpha(v,s)}=\deltap/\deltam$, which shows that most paths
starting at level $1-\deltap$ at time $s$ will reach the lower level
$1-\deltam$ close to that time. 

The same argument as in Proposition~\ref{prop_fps}, applied for a
different initial condition, shows that $\psiup(u,v)$ can be written as
$\sigma^{-1}\ceeup(u,v)\e^{-\rhoup(u,v)^2/2\sigma^2}$, with a rate 
\begin{equation}
\label{fpK2}
\rhoup(u,v)^2 = 
\frac{\bigbrak{(2-\deltap)-(2-\deltam)\e^{-\alpha(u,v)}}^2}{v_-(u,v)},
\end{equation}
which is bounded below by $(\deltam-\deltap)^2/\bv$. It thus follows
from the definition~\eqref{re5b} of $K$ that $K(u,s)\leqs
(\const/\sigma)\e^{-(\deltam-\deltap)^2/2\overline v\sigma^2}$. 

By differentiating~\eqref{re10}, we find  
\begin{equation}
\label{fpK10}
q(t,s) = p_1(t,s) + \sum_{n=1}^\infty \int_s^t p_1(t,u)K_n(u,s)\6u.
\end{equation} 
(Note that the sum is converging uniformly on compact sets.)
The smallness of $K$ is not yet sufficient to ensure the smallness of the
sum, relatively to $p_1(t,s)$. The following result provides a sufficient
bound. 

\begin{prop}
\label{prop_fpK}
Let $\overline p_1(t,s) \defby \frac1\sigma (c_+(t,s)\vee1)
\e^{-\rho_+(t,s)^2/2\sigma^2} \geqs p_1(t,s)$.\footnote{We write
$a\vee b$ to denote the maximum of two real numbers $a$ and $b$.}\ %
Assume that $\vhatper_+(u) \leqs 2 v^\star(u)(1-\Delta)$ for all~$u$. 
Then, the relation
\begin{equation}
\label{fpK11}
\int_s^t \overline p_1(t,u) K_n(u,s)\6u 
\leqs \biggbrak{\const (t-s)\frac1{\sigma^2}\e^{-\const\Delta_0^2/\sigma^2}}^n 
\overline p_1(t,s)
\end{equation}
holds for all $n\geqs1$, with $\Delta_0 =
\Delta/(1-\Delta)\wedge(\deltam-\deltap)/\deltap\wedge1$. 
\end{prop}
\begin{proof}
The main idea consists in comparing $\psiup(u,v)$ to the density
$\psiup^+(u,v)$ of the first-passage time of $y^+_t(v,1-\deltam)$ at
$1-\deltap$ (recall that $\psiup(u,v)$ relates to the process
$y^-_t(v,1-\deltam)$). The exponential rate $\rhoup^+(u,v)$ associated with
$\psiup^+(u,v)$ is given by 
\begin{equation}
\label{fpK12:1}
\rhoup^+(u,v)^2 = 
\frac{\bigbrak{\deltam-\deltap\e^{-\alpha(u,v)}}^2}{\vhat_+(u,v)}.
\end{equation}
We claim that 
\begin{equation}
\label{fpK12:2}
\rhoup(u,v)^2 \geqs \rhoup^+(u,v)^2
\end{equation}
for all $u\geqs v$, as a consequence of Hypothesis~H4. 
Denote $\e^{-\alpha(u,v)}$ by $\xi$ and $\underline v/\bv$ by $\kappa^2$.
On the one hand, it follows from the definitions~\eqref{fps3}
and~\eqref{fpu4} of $v_-$ and $\vhat_+$ that 
\begin{equation}
\label{fpK12:3}
v_-(u,v) \leqs \int_v^u g(w)^2 \6w \leqs \frac1{\xi^2} \vhat_+(u,v). 
\end{equation}
On the other hand, since $g(w)^2=2a(w)v^\star(w)$ and
$v^\star(w)\in[\underline v,\bv]$, we can write 
\begin{equation}
\label{fpK12:4}
\frac1{\bv}v_-(u,v) \leqs \int_v^u 2a(w)\e^{-2\alpha(u,w)} \6w
= 1 - \xi^2 = \int_v^u 2a(w)\e^{-2\alpha(w,v)} \6w
\leqs \frac1{\underline v} \vhat_+(u,v). 
\end{equation}
We can thus conclude that 
\begin{equation}
\label{fpK12:5}
\frac{\rhoup(u,v)}{\rhoup^+(u,v)} \geqs 
\frac{(2-\deltap)-(2-\deltam)\xi}{\deltam - \deltap\xi} 
\bigbrak{\xi\vee\kappa}. 
\end{equation}
Using the monotonicity of the ratio in the preceding estimate, we see 
that Relation~\eqref{fpK12:2} is satisfied provided
$\brak{(2-\deltap)-(2-\deltam)\xi}\xi \geqs \deltam - \deltap\xi$ for all
$\xi\geqs\kappa$, which easily follows from
Hypothesis~H4. 

In order to prove~\eqref{fpK11} for $n=1$, we first consider the rate
$\tilde\rho_+(t,v)$ associated with $\probin{v,1-\deltam}{\tilde\tau_+
\leqs t}$, compare~\eqref{fpu18:4}. Proceeding as in
Proposition~\ref{prop_fpu}, it is straightforward to show that 
\begin{align}
\nonumber
\tilde\rho_+(t,v)^2 &\leqs \rho_+(t,u)^2 + \rhoup^+(u,v)^2 \\
&\leqs \rho_+(t,u)^2 + \rhoup(u,v)^2
\label{fpK12:7}
\end{align}
for all $u\in[v,t]$. This is actually a consequence of the Markov property
and the fact that $y^+_t(v,1-\delta_2)$ has to cross the level $1-\delta_1$
before reaching $1+$. Equation~\eqref{fpK12:7} implies that 
\begin{align}
\nonumber
\int_s^t \overline p_1(t,u) & K(u,s)\6u \\
\nonumber
&\leqs \int_s^t\int_v^t \frac1{\sigma^2} \bigbrak{c_+(t,u)\vee1}\ceeup(u,v)
\e^{-\brak{\rho_+(t,u)^2+\rhoup(u,v)^2}/2\sigma^2} \6u \mskip6mu \psidown(v,s)\6v
\\
\nonumber
&\leqs \const \frac{t-s}{\sigma^2}
\int_s^t\e^{-\tilde\rho_+(t,v)^2/2\sigma^2}\psidown(v,s)\6v \\
&\leqs \const \frac{t-s}{\sigma^3} \e^{-\const\Delta_0^2/\sigma^2} 
\e^{-\rho_+(t,s)^2/2\sigma^2},
\label{fpK12:8}
\end{align}
where the last line is obtained by the same argument as in
Proposition~\ref{prop_fpu}, compare~\eqref{fpu18:4b}. This
proves~\eqref{fpK11} for $n=1$, and for larger $n$ the result follows easily
by induction. 
\end{proof}

As a direct consequence of Proposition~\ref{prop_fpK}, whenever $t-s < 
\sigma^2 \e^{\const\Delta_0^2/\sigma^2}$, we can bound the sum
in~\eqref{fpK10} by a geometric series, and conclude that   
\begin{equation}
\label{fpK13}
\frac1\sigma c_+(t,s) \e^{-\rho_+(t,s)^2/2\sigma^2} 
\leqs q(t,s) 
\leqs \frac1\sigma \bar c_+(t,s) \e^{-\rho_+(t,s)^2/2\sigma^2},
\end{equation}
where $\bar c_+(t,s) = c_+(t,s) +
\bigOrder{(t-s)\sigma^{-2}(1 + c_+(t,s))\e^{-\const\Delta_0^2/\sigma^2}}$. 


\section{Properties of the exit law}
\label{sec_xl}


The first-passage law of the process $y_t$ at $+1$ can be expressed as a
function of the first-passage density $\psi_-(s,0)$ to $1-\deltap$ and
$q(t,s)$ via~\eqref{re1} as 
\begin{equation}
p_+(t) = \dpar{}{t} \probin{0,-1}{\tau_+ \leqs t} 
=\int_0^t q(t,s)\psi_-(s,0)\6s.
\label{xl1}
\end{equation}
From now on, we assume that
$t\leqs\sigma^3\e^{\beta\Delta_0^2/2\sigma^2}$ for some sufficiently
small constant $\beta>0$, so that the error terms
in~\eqref{fps10} and~\eqref{fpu17} are at most of order $h(\sigma) =
\sigma \e^{-\beta\Delta_0^2/2\sigma^2}$, and the one in~\eqref{fpK13}
is at most of order $h(\sigma)(1+c_+(t,s))$.


\subsection{Estimating the integral by a sum}
\label{ssec_xll}

\begin{figure}
\centerline{\includegraphics*[clip=true,width=125mm]{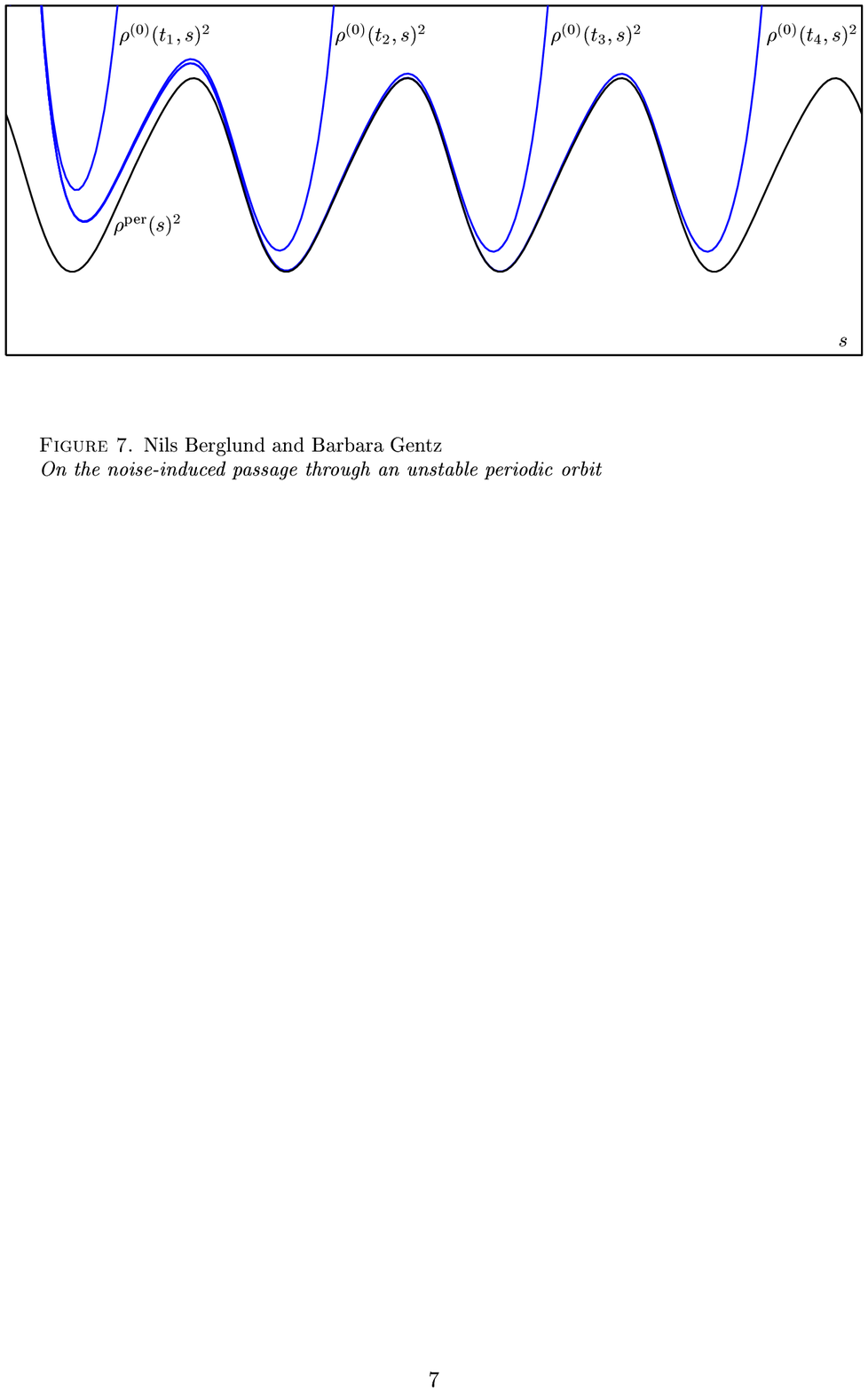}}
 \captionspace
 \caption[]
 {The functions $\rhoper(s)^2$ and $\rho^{(0)}(t_i,s)^2$ for four different
 final times $t_1 < t_2 < t_3 < t_4$.}
\label{fig7}
\end{figure}

In virtue of Proposition~\ref{prop_fps}, 
Proposition~\ref{prop_fpu}, and~\eqref{fpK13}, the integral in~\eqref{xl1}
can be written as 
\begin{equation}
p_+(t) 
= \frac1{\sigma^2} \int_0^t \bigbrak{
  c_+(t,s)+\bigOrder{h(\sigma)[1+c_+(t,s)]}} c_-(s,0) 
\e^{-[\rho_+(t,s)^2+\rho_-(s,0)^2]/2\sigma^2}\6s.
\label{xll1}
\end{equation}
The exponent in~\eqref{xll1} is of the form 
\begin{align}
\nonumber
\rho^{(0)}(t,s)^2 &:= \rho_+(t,s)^2+\rho_-(s,0)^2 \\
\nonumber
&=\frac{\deltap^2}{\vhat_+(t,s)} + 
\frac{(2-\deltap)^2}{v_-(s,0)} \\
\label{xll2a}
&=\frac{\deltap^2}{\vhatper_+(s)-\e^{-2\alpha(t,s)}\vhatper_+(t)} + 
\frac{(2-\deltap)^2}{\vper_-(s)-\e^{-2\alpha(s)}\vper_-(0)},
\end{align}
which we write for $0\ll s\ll t$ as
\begin{equation}
\label{xll2b}
\rho^{(0)}(t,s)^2 =\rhoper(s)^2 + \Order{\e^{-2\alpha(t,s)}} +
\Order{\e^{-2\alpha(s)}}
\end{equation}
with a periodic part 
\begin{equation}
\label{xll3}
\rhoper(s)^2 = \frac{\deltap^2}{\vhatper_+(s)} 
+ \frac{(2-\deltap)^2}{\vper_-(s)},
\end{equation}
cf.~\figref{fig7}. In Hypothesis H5, we assumed, for simplicity,  that
$\rhoper(s)$ has a unique minimum in the interval $[0,T)$, at some
$s^\star\in[0,T)$ satisfying 
\begin{equation}
\label{xll4}
\rhoper(s^\star) = \inf_{s\in[0,T)} \rhoper(s). 
\end{equation}
When the time interval $[0,t]$ includes many periods, $\rho^{(0)}(t,s)$
will have minima of comparable depths near all times $0 \ll s^\star+kT \ll
t$, see~\figref{fig7}. Other minima, which may exist for $s$ near $0$ and $t$
are much shallower, and thus contribute less to the integral~\eqref{xll1}.
The integral will be of the order 
$(t/T)\sigma \e^{-\rhoper(s^\star)^2/2\sigma^2}$. 
On  the other hand, when $t$ is not large enough, even the deepest minimum of
$\rho^{(0)}(t,s)^2$ will be substantially larger than $\rhoper(s^\star)^2$,
leading to a value of $p_+(t)$ which is orders of magnitude smaller
than 
$(t/T)\sigma \e^{-\rhoper(s^\star)^2/2\sigma^2}$. 
The system is then still in its initial transient regime. 

The transition between these regimes occurs when
$\alpha(t)\simeq2\abs{\log\sigma}$. We first show
that for times such that $\alpha(t)$ is smaller than
$2\abs{\log\sigma}$, the density 
$p_+(t)$ is much smaller than its \lq\lq asymptotic\rq\rq\ value
$(t/\sigma T)\e^{-\rhoper(s^\star)^2/2\sigma^2}$.

\begin{prop}
\label{prop_xll1}
Assume that $\alpha(t) \leqs 2\nu\abs{\log\sigma}$ for some $\nu<1$. Then
\begin{equation}
\label{xll4a}
p_+(t) \leqs \const\frac{1}{\sigma^2}
\e^{-L/\sigma^{2(1-\nu)}} \e^{-\rhoper(s^\star)^2/2\sigma^2},
\end{equation}  
where $L>0$ is a constant independent of $\sigma$, $t$ and $\nu$. 
\end{prop}
\begin{proof}
From~\eqref{xll2a} we obtain 
\begin{equation}
\label{xll4b}
\rho^{(0)}(t,s)^2 \geqs \rhoper(s)^2 + 
\deltap^2\frac{\vhatper_+(t)}{\vhatper_+(s)^2}\e^{-2\alpha(t,s)} + 
(2-\deltap)^2\frac{\vper_-(0)}{\vper_-(s)^2}\e^{-2\alpha(s)}.
\end{equation}
The periodic functions $\vhatper_+(s)$ and $\vper_-(s)$ being bounded above
by $\bv$, we have 
\begin{equation}
\label{xll4c}
\rho^{(0)}(t,s)^2 \geqs \rhoper(s^\star)^2 + \beta_1^2 \e^{-2\alpha(t,s)} +
\beta_2^2\e^{-2\alpha(s)}, 
\end{equation}
where $\beta_1=\deltap\vhatper_+(t)^{1/2}/\bar{v}$ and 
$\beta_2=(2-\deltap)\vper_-(0)^{1/2}/\bar{v}$. The right-hand side
of~\eqref{xll4c} reaches its minimum when~$s$ satisfies 
$\beta_2\e^{-2\alpha(s)}=\beta_1\e^{-\alpha(t)}$, and has value 
\begin{equation}
\label{xll4d}
\rhoper(s^\star)^2 + 2\beta_1\beta_2\e^{-2\alpha(t)} 
\geqs\rhoper(s^\star)^2 + 2\beta_1\beta_2\sigma^{2\nu}.
\end{equation}
We can now estimate the integral~\eqref{xll1} by splitting it at times
$t_1=1$ and $t_2=t-1$. (If $t\leqs 2$, the argument is even simpler.)
The integral over $[t_1,t_2]$ can be bounded by using~\eqref{xll4d}, while the
integral on $[0,t_1]$ is small because 
for $\sigma$ small enough, $c_-(s,0)\e^{-\rho_-(s,0)^2/2\sigma^2}$ reaches  
its maximum at $t_1$ and the remaining factor in the integrand is bounded.
On $[t_2,t]$ the situation is similar. 
\end{proof}

We now turn to the case where $t>2\abs{\log\sigma}/\lambda$. Assume that
$t\in[nT,(n+1)T)$. If we can show that $\rho^{(0)}(t,s)$ has exactly one
minimum $s_k$ in $[kT,(k+1)T)$ for $1\ll k\ll n$, and that this minimum is
quadratic, the Laplace method will allow us to approximate the
integral~\eqref{xll1} by 
\begin{equation}
\label{xll5}
\frac1\sigma \biggpar{\sum_{k} 
\sqrt{\frac{4\pi}{\sdpar{}{ss}(\rho^{(0)}(t,s_k)^2)}}\mskip1.5mu
c_+(t,s_k) c_-(s_k,0) + \text{\it error term}} \e^{-\rho^{(0)}(t,s_k)^2/2\sigma^2}. 
\end{equation}
Making this argument precise, we obtain the following result. 

\begin{prop}
\label{prop_xll2}
Assume that $t\in[nT,(n+1)T)$ with $n\lambda T \geqs 2\abs{\log\sigma}$ and
$n\geqs4$. Then 
\begin{equation}
\label{xll5b}
p_+(t) = \frac1\sigma C(s^\star) \frac{g(t)^2}{\vhatper_+(t)}
S(n,\sigma,t) \bigbrak{1+\Order{\sigma}}
\e^{-\rhoper(s^\star)^2/2\sigma^2},
\end{equation}
where $C(s^\star)$ is a constant given by 
\begin{equation}
\label{xll15d}
C(s^\star) = \frac{2(2-\deltap)}{\deltap}
\frac{g(s^\star)^2}{\sqrt{\pi\sdpar{}{ss}(\rhoper(s^\star)^2)}}
\frac{\vhatper_+(s^\star)^{1/2}}{\vper_-(s^\star)^{3/2}}
\biggbrak{1-\frac{\vper_-(s^\star)}{2v^\star(s^\star)}}, 
\end{equation}
and $S(n,\sigma,t)$ is the sum
\begin{equation}
S(n,\sigma,t) = \frac{\gamma(t)}2
\sum_{k=1}^n \exp\biggset{-2(n-k)\lambda T - \frac1{2\sigma^2}
\Bigbrak{\gamma_0\e^{-2k\lambda T} + \gamma(t)\e^{-2(n-k)\lambda T}}}, 
\label{xll6}
\end{equation}
with  
\begin{align}
\nonumber
\gamma_0 &= (2-\deltap)^2 \e^{-2\alpha(s^\star)}
\frac{\vper_-(0)}{\vper_-(s^\star)^2},\\
\label{xll7}
\gamma(t) &= \deltap^2 \e^{-2\alpha(t,s^\star+nT)}
\frac{\vhatper_+(t)}{\vhatper_+(s^\star)^2}.
\end{align}
\end{prop}
\begin{proof}
We split the integral~\eqref{xll1} at times $t_1=k_1 T$ and $t_2=(n-k_1)T$,
where $k_1<n/2$ will be chosen in such a way that for $t_1\leqs s\leqs t_2$,
$\rho^{(0)}(t,s)^2$ has a minimum close to $s^\star+kT$ on each interval
$I_k=[kT,(k+1)T)$, while the contributions of the integrals over $[0,t_1]$ and
$[t_2,t]$ are negligible. 

Take $k_1$ of the form $k_1=\intpartplus{\nu\abs{\log\sigma}/2\lambda T}
\vee 2$, with a parameter $\nu$ yet to be chosen. For $t_1\leqs s\leqs
t_2$, one has  $\e^{-2\alpha(s)}\leqs\sigma^\nu$ and
$\e^{-2\alpha(t,s)}\leqs\sigma^\nu$.  We first show that
$\rho^{(0)}(t,s)^2$ has a quadratic minimum in each $I_k$, $k_1\leqs
k<n-k_1$. For this purpose, we write 
\begin{equation}
\label{xll20:1}
\rho^{(0)}(t,s)^2 =
\frac{\deltap^2}{\vhatper_+(s)-\e^{2\alpha(s,kT)}\gamma_1} + 
\frac{(2-\deltap)^2}{\vper_-(s)-\e^{-2\alpha(s,kT)}\gamma_2},
\end{equation}
where $\gamma_1 = \e^{-2(n-k)\lambda T}\e^{-2\alpha(t,nT)}\vhatper_+(t)$ and
$\gamma_2 = \e^{-2k\lambda T}\vper_-(0)$ are at most of order
$\sigma^\nu$. Were $\gamma_1=\gamma_2=0$, \eqref{xll20:1} would
reduce to $\rhoper(s)^2$, which has 
a unique minimum in $I_k$ at $s^\star+kT$, the latter being
quadratic. Hence, the implicit-function theorem applied to
$\sdpar{}{s}\rho^{(0)}(t,s)^2$ shows that for sufficiently small
$\gamma_1$ and $\gamma_2$, $\rho^{(0)}(t,s)^2$ has a unique minimum in
$I_k$ at a time $s_k=s^\star+kT+\Order{\sigma^\nu}$. 

Expanding~\eqref{xll20:1} into powers of $\gamma_1$ and $\gamma_2$ shows
that 
\begin{equation}
\label{xll20:2}
\rho^{(0)}(t,s_k)^2 = \rhoper(s^\star)^2 +
\gamma(t)\e^{-2(n-k)\lambda T}
+ \gamma_0\e^{-2k\lambda T}  + \Order{\sigma^{2\nu}},  
\end{equation}
for the coefficients $\gamma_0$ and $\gamma(t)$ given in~\eqref{xll7}. 
Evaluating the integral~\eqref{xll1} restricted to the interval $[t_1,t_2]$
by the Laplace method yields 
\begin{equation}
\label{xll20:3}
\begin{split}
\frac1\sigma \sum_{k=k_1}^{n-k_1-1}
&{}\sqrt{\frac{4\pi}{\sdpar{}{ss}(\rhoper(s_k)^2)}} 
\bigbrak{
  c_+(t,s_k) + \bigOrder{h(\sigma)[1+c_+(t,s_k)]}} 
c_-(s_k,0)\bigbrak{1+\Order{\sigma^2}} \\
&{}\times \e^{-\rho^{(0)}(t,s_k)^2/2\sigma^2},
\end{split}
\end{equation}
where $\sdpar{}{ss}\rhoper(s_k) = \sdpar{}{ss}\rhoper(s^\star) +
\Order{\sigma^\nu}$ and 
\begin{align}
\label{xll20:4a}
c_+(t,s_k) &= \frac{\deltap}{\sqrt{2\pi}} 
\frac{g(t)^2\e^{-2\alpha(t,s^\star+nT)}\e^{-2(n-k)\lambda T}}
{\vhatper_+(s^\star)^{3/2}}\bigbrak{1+\Order{\sigma^\nu}} -
\Order{h(\sigma)},\\ 
\label{xll20:4b}
c_-(s_k,0) &= \frac{2-\deltap}{\sqrt{2\pi}} 
\biggbrak{\frac1{\vper_-(s^\star)}-\frac1{2v^\star(s^\star)}}
\frac{g(s^\star)^2}{\vper_-(s^\star)^{1/2}}
\bigbrak{1+\Order{\sigma^\nu+h(\sigma)}}.
\end{align}
The error term $\Order{\sigma^{2\nu}}$ in the exponent~\eqref{xll20:2}
yields an error term $\Order{\sigma^{2(\nu-1)}}$ in the sum~\eqref{xll20:3}.
For this reason, we choose $\nu=3/2$, so that all error terms are of order
$\sigma$ at most. This choice of $\nu$ is always possible, since $2k_1 \leqs
3\abs{\log\sigma}/(2\lambda T)+2$ is smaller than $n$ by our condition on
$\abs{\log\sigma}$. Note that the additive error term
$\Order{h(\sigma)}$ in~\eqref{xll20:4a} can be incorporated into the
sum $S(n,\sigma,t)$, as $nh(\sigma) \leqs \Order{th(\sigma)} \leqs
\Order{\sigma^4} \leqs \Order{\sigma^2} S(n,\sigma,t)$, where the last
inequality will be proved in the next section. 

We now turn to computing a bound for the integral~\eqref{xll1} restricted
to the interval $[0,t_1]$. We first consider the case 
$\abs{\log\sigma}\geqs 5\lambda T$, in which $k_1\geqs4$ and
$\e^{-2\alpha(t_1)}=\e^{-2k_1\lambda T}\geqs \sigma^{8/5}$. 
With the help of~\eqref{xll4c}, the integral can be estimated by
$\exp\{-\rhoper(s^\star)^2/2\sigma^2\}$ times its maximal value which itself
is bounded by
\begin{equation}
\frac\const{t_1^{1/2}\sigma^2} 
\exp\biggset{-\frac{\beta_2^2}{2\sigma^2}\e^{-2\alpha(t_1)}}
=
\biggOrder{\frac1{\sigma^{2}\abs{\log\sigma}^{1/2}}
\exp\Bigset{-\frac{\beta_2^2}{2\sigma^{2/5}}}}.
\label{xll20:5}
\end{equation}
Hence the first part of the integral is exponentially small compared
to the integral itself.  A similar
estimate shows that the integral over $[t_2,t]$ is small.  

In the case $\abs{\log\sigma}<5\lambda T$, the location of the minima has
to be estimated with more care. In this case, $\e^{-2\lambda
T}\leqs\sigma^{2/5}$, so that the same argument as before shows that
$\rho^{(0)}(t,s)^2$ has a minimum in each of the intervals $I_1$,
$I_2$ and $I_3$, and possibly 
also in $I_0$. An examination of $\sdpar{}s \rho^{(0)}(t,s)^2$ shows that
these minima are actually located in
$s_k=s^\star+kT+\Order{\e^{-2\alpha(s^\star+kT)}}$, $k=0,1,2,3$. We set
$\kappa=\alpha(s^\star)/\lambda T$, $0\leqs\kappa<1$, and distinguish
between two cases: 
\begin{itemiz}
\item	If $4(k+\kappa)\lambda T\geqs3\abs{\log\sigma}$, then 
$\rho^{(0)}(t,s_k)^2=\rhoper(s^\star)^2+\gamma_0\e^{-2k\lambda T} +
\Order{\sigma^3}$, and the integral between $0$ and
$s_k$ is exponentially small.
\item	If $4(k+\kappa)\lambda T<3\abs{\log\sigma}$, then
$\e^{-2\alpha(s_k)}\geqs\sigma^{3/2}$, and the error term only leads to an
error of order $\sigma$ in the prefactor. 
\end{itemiz}
The intervals $I_k$, $n-4\leqs k\leqs n-1$, are treated in a similar way. 

The same type of arguments shows that the error made by extending the sum
in~\eqref{xll20:3} to all $k$ between $1$ and $n$ is also negligible. 
\end{proof}


\subsection{Properties of the sum}
\label{ssec_xls}

To complete the proof of Theorem~\ref{thm_mrr}, we will show that the sum
$S(n,\sigma,t)$, defined in~\eqref{xll6}, is close to a periodic function
of $\abs{\log\sigma}$.  In this section, we always assume that
$t\in[nT,(n+1)T)$ and that $n\lambda T \geqs 2\abs{\log\sigma}$ as well as
$n\geqs4$. To highlight the  $\sigma$-dependence of $S$, it is convenient
to set $\sigma = \e^{-\eta}$ and $S(n,\sigma,t)=\sigma^2
\Stilde(n,\eta,t)$. Further introducing the notations 
\begin{equation}
\label{xls0}
\theta_0 = -\frac12 \log\gamma_0, \qquad
\bar\theta(t) = -\frac12 \log\gamma(t),
\end{equation}
and changing the summation index from $k$ to $\l=n-k$, allows to write the
sum in compact form as 
\begin{equation}
\label{xls1}
\Stilde(n,\eta,t) = \sum_{\l=0}^{n-1} A\bigpar{\l\lambda T-\eta+\bar\theta(t)}
B\bigpar{(n-\l)\lambda T-\eta+\theta_0},
\end{equation}
with 
\begin{align}
\label{xls2}
A(x) &= \frac12\exp\biggset{-2x - \frac12 \e^{-2x}}, \\
\label{xls3}
B(x) &= \exp\biggset{-\frac12 \e^{-2x}}. 
\end{align}
The function $A(x)$ decays like $\e^{-2x}$ as $x\to +\infty$,
and like $\exp\set{-\frac12\e^{2\abs{x}}}$ as $x\to -\infty$. 
It reaches its maximal value $\e^{-1}$ when $\e^{-2x}=2$. 
The function $B(x)$ is monotonously increasing. It decays like
$\exp\set{-\frac12\e^{2\abs{x}}}$ as $x\to -\infty$ and
approaches $1$ as $x\to+\infty$.

Since $\gamma(t)$ involves the location of $t\in[nT,(n+1)T)$ relatively to
$s^\star+nT$ (see~\eqref{xll7}), $\bar\theta(t)$ is a (right-continuous)
periodic saw-tooth function, making a jump of $-\lambda T$ at each integer
multiple of $T$. We can thus write it as 
\begin{equation}
\label{xls3a}
\bar\theta(t) = \theta(t) - \lambda T \biggintpart{\frac{t}{T}}, 
\end{equation}
where $\intpart{x}$ denotes the largest integer smaller than or equal to
$x$, and $\theta(t)$ is a continuous function given by 
\begin{equation}
\label{xls3b}
\theta(t) = \alpha(t,s^\star) - \frac12 \log\vhatper_+(t) 
- \log\frac{2-\deltap}{\vhatper_+(s^\star)}. 
\end{equation}
Note that $\theta(t+T) = \theta(t) + \lambda T$. 

It is easy to show that $\Stilde(n,\eta+\lambda T,t) -
\Stilde(n-2,\eta,t-2T)$ is exponentially small (just shift the index
of summation and show that the boundary terms are exponentially
small). The following result gives a more precise characterization of
$\Stilde(n,\eta,t)$ by showing that it is actually close to a periodic
function of $\eta-\theta(t)$. 

\begin{prop}
\label{prop_xls1}
Assume that $n\lambda T \geqs 2\eta$. Then 
\begin{equation}
\label{xls7}
\Stilde(n,\eta,t) = \Shat(\eta,t) \bigbrak{1 + \Order{\sigma^\mu}}, 
\end{equation}
where $\mu = \mu(n,\eta) = (n\lambda T-2\eta)/\eta$ and  
\begin{equation}
\label{xls6}
\Shat(\eta,t) = 
\sum_{\l=-\infty}^{\infty} A\bigpar{\l\lambda T-\eta+\theta(t)}. 
\end{equation}
\end{prop}
\begin{proof}
We split the sum $\Stilde(n,\eta,t)$ into two parts. For
$0\leqs\l\leqs n/2$, we use the fact that 
\begin{equation}
\label{xls10:2}
1\geqs B\bigpar{(n-\l)\lambda T-\eta+\theta_0} 
\geqs B\bigpar{\tfrac n2\lambda T-\eta+\theta_0}
= \exp\Bigset{-\frac{\gamma_0}2\e^{2\eta-n\lambda T}}
= 1 - \Order{\sigma^\mu}. 
\end{equation}
Hence replacing $B$ by $1$ in the sum for these values of $\l$ only
yields a multiplicative error of order $1 - \Order{\sigma^\mu}$. 
For $\l>n/2$, it is obvious that 
\begin{equation}
\label{xls10:3}
A\bigpar{\l\lambda T-\eta+\bar\theta(t)}
\leqs \frac{\gamma(t)}2 \e^{-2\l\lambda T+2\eta}.
\end{equation}
Bounding $B$ by $1$ allows to bound the sum over
$\l = \intpartplus{n/2},\dots,n$ by the geometric series   
\begin{equation}
\label{xls10:4}
\sum_{\l=\intpartplus{n/2}}^n A\bigpar{\l\lambda T-\eta+\bar\theta(t)} 
\leqs \frac{\gamma(t)}2\sum_{\l=\intpartplus{n/2}}^n \e^{-2\l\lambda T+2\eta}
= \Order{\sigma^\mu}. 
\end{equation}
Thus the main contribution to $\Stilde$ stems from
$\l\in\{1,\dots,\intpart{n/2}\}$, and so does the main contribution to
$\sum_{\l=0}^{\infty} A\bigpar{\l\lambda T-\eta+\theta(t)}$.
It remains to check that the contribution of negative $\l$ to $\Shat$
is small. Comparing that sum with an
integral shows that it is in fact of the order
$\e^{-\gamma(t)/2\sigma^2}$. Finally, replacing $\bar\theta(t)$ by
$\theta(t)$ only results in a shift of the summation index. 
\end{proof}

The function $\Shat(\eta,t)$ is clearly periodic in $\eta-\theta(t)$
with period $\lambda T$. Let us thus write $\Shat(\eta,t) =
P((\eta-\theta(t))/\lambda T)$, where $P(x)>0$ is periodic with period $1$. 
It remains to compute the Fourier series of $P(x)$. 

\begin{prop}
\label{prop_xls2}
The periodic function $P(x)$ admits the Fourier series 
\begin{equation}
\label{xls11}
P(x) = 
\sum_{q=-\infty}^{\infty} \Phat(q) \e^{2\pi\icx qx}, 
\end{equation}
where the $q$th Fourier coefficient is given in terms of the Euler Gamma
function by 
\begin{equation}
\label{xls12}
\Phat(q) = \frac1{2\lambda T}
\frac{1}{2^{\pi\icx q/\lambda T}} 
\Gamma\biggpar{1 - \frac{\pi\icx q}{\lambda T}}. 
\end{equation}
\end{prop}
\begin{proof}
We have 
\begin{align}
\nonumber
\Phat(q) &=  \int_0^1 
P(x)\e^{-2\pi\icx qx} \6x \\
\nonumber
&= \sum_{\l=-\infty}^\infty \int_0^1
A\bigpar{(\l-x)\lambda T}\e^{-2\pi\icx qx} \6x \\
&= \int_{-\infty}^\infty A(\lambda Tx) \e^{2\pi\icx qx}\6x.
\label{xls13}
\end{align}
Replacing $A$ by its definition and using the change of variable 
$z=\e^{-2\lambda Tx}/2$ yields the result. 
\end{proof}



\appendix
 
\makeatletter
\renewcommand\theequation{A.\@arabic\c@equation}
\renewcommand\thetheorem{A.\@arabic\c@theorem}
\renewcommand\thesection{Appendix}
\makeatother
 
\setcounter{equation}{0}
\setcounter{theorem}{0}    

\section{}
\label{appendix}


Let $v(t)$ be continuously differentiable on $[0,\infty)$ and satisfy
$v(0)=0$ and $\vprime(t)\geqs v_0>0$ for all $t\geqs 0$. As before, we
denote by $W_t$ a standard Brownian motion. Consider the Gaussian
process
\begin{equation}
\label{app0}
z_t = \sigma \int_0^t \sqrt{\vprime(s)}\6W_s, 
\end{equation}
whose variance is $\sigma^2 v(t)$. We will consider $z_t$ as a Markov
process and introduce the notation $\probin{s,x}{z_t \in \cdot} =
\pcond{z_t \in \cdot\,}{z_s = x}$, $t>s$, for its transition
probabilities. Their densities are given by
\begin{align}
\nonumber
y\mapsto f(t,y \vert s,x) 
&\defby \dpar{}y \probin{s,x}{z_t \leqs y} \\
&=\frac1\sigma\frac1{\sqrt{2\pi v(t,s)}} \e^{-(y-x)^2/2\sigma^2v(t,s)}
\label{app0a}
\end{align}
where $v(t,s)=v(t)-v(s)$. 

Let $d(t)$ be continuously differentiable on $[0,\infty)$ and satisfy
$d(0)>0$. The object of the \lq\lq level-crossing problem\rq\rq\ is to
determine the density $\psi(t)$ of the first-passage time
$\tau=\inf\setsuch{s>0}{z_s>d(s)}$ (which we will call \lq\lq
first-passage density of $z_t$ to $d(t)$\rq\rq). This problem has for
instance been studied in~\cite{Durbin85,RicciSacer84,Durbin,SacerTomas96}. 
The aim of this appendix is to establish expressions for $\psi(t)$ useful
in our particular situation.  

As D.~Williams in the appendix to~\cite{Durbin}, we will use two
integral equations satisfied by $\psi(t)$. Let $F(t)=f(t,d(t)\vert
0,0)$ denote the value of the density of $z_t$ at $d(t)$ and let
$F(t\vert s) = f(t,d(t)\vert s,d(s))$ denote the transition density at
$y=d(t)$ for paths starting at time $s$ in $x=d(s)$. The Markov
property enables us to write  
\begin{equation}
\label{app0c}
F(t) = \int_0^t F(t\vert s) \psi(s) \6s.
\end{equation}
The second integral equation satisfied by $\psi(t)$ is 
\begin{equation}
\label{app1}
\psi(t) = b_0(t) F(t)
- \int_0^t \tilde b(t,s) F(t\vert s) \psi(s)\6s,
\end{equation}
where
\begin{align}
\label{app2a}
b_0(t) &= \vprime(t) 
\biggbrak{\frac{d(t)}{v(t)} - \frac{\dprime(t)}{\vprime(t)}}
= - v(t) \dpar{}{t}\biggpar{\frac{d(t)}{v(t)}}, \\
\label{app2b}
\tilde b(t,s) &= \vprime(t)
\biggbrak{\frac{d(t,s)}{v(t,s)} - \frac{\dprime(t)}{\vprime(t)}}
= - v(t,s) \dpar{}{t}\biggpar{\frac{d(t,s)}{v(t,s)}},
\end{align}
with $d(t,s)=d(t)-d(s)$. In the particular case of a standard Brownian
motion, that is for $\sigma=1$ and $v(t)\equiv t$,
Equation~\eqref{app1} has been established in~\cite[Appendix by
  D.~Williams]{Durbin}. The general case is easily obtained from the fact
that $z_t=\sigma W_{v(t)}$ in distribution.

Equation~\eqref{app1} suggests that the first-passage density can be written
in the form 
\begin{equation}
\label{app2c}
\psi(t) = \frac1\sigma c(t) \e^{-d(t)^2/2\sigma^2 v(t)}, 
\end{equation}
where $c(t)$ is a subexponential prefactor. In fact, the following bound on
$c(t)$ follows immediately from~\eqref{app0c} and~\eqref{app1}. 

\begin{lemma}
\label{lem_appendix0}
Let $c_0(t) = b_0(t)/\sqrt{2\pi v(t)}$. Then~\eqref{app2c} holds with
\begin{equation}
\label{app2d}
\bigabs{c(t)-c_0(t)} \leqs \frac1{\sqrt{2\pi v(t)}} \sup_{0\leqs s\leqs
t}\bigabs{\tilde b(t,s)}. 
\end{equation}
\end{lemma}

\begin{remark}
\label{rem_appendix}
Note that Lemma~\ref{lem_appendix0} does not require $\vprime(t)$ to
be bounded away from zero as $t$ varies.
\end{remark}

If, for instance, $v(t)$ and $d(t)$ are twice continuously differentiable,
then $s\mapsto\tilde b(t,s)$ is easily seen to be bounded, and thus
$c(t)$ behaves like $v(t)^{-3/2}$ near $t=0$.

In~\cite{Durbin}, an expansion of $c(t)$ is constructed, which is shown to
converge for all times $t$, under a convexity assumption on the boundary
$d(t)$. Taking advantage of the fact that $\sigma$ is a small
parameter, we can control the convergence of this expansion under a
milder assumption on $d(t)$, on a finite, but exponentially long time
interval. Writing $\tilde c(t,s) = \tilde b(t,s) / \sqrt{2\pi
v(t,s)}$, we see from~\eqref{app1} that $c(t)$ must be a fixed point
of the operator 
\begin{equation}
\label{app3}
(\cT c)(t) = c_0(t) - \frac1\sigma  
\int_0^t \tilde c(t,s)c(s) \e^{-r(t,s)/2\sigma^2}\6s,
\end{equation}
where
\begin{equation}
\label{app4}
r(t,s) 
= \frac{d(s)^2}{v(s)} - \frac{d(t)^2}{v(t)} + \frac{d(t,s)^2}{v(t,s)} 
= \frac{v(t)v(s)}{v(t,s)} \biggbrak{\frac{d(s)}{v(s)} - \frac{d(t)}{v(t)}}^2.
\end{equation}

\begin{remark}   \label{rem2_appendix}
The exponent $r(t,s)$ is nonnegative and vanishes for $s=t$. If $r(t,s)$
does not vanish anywhere else, then the main contribution to the  integral
in~\eqref{app3} comes from $s$ close to $t$. In the generic case $\sdpar
rs(t,t)\neq0$, the integral is at most of order $\sigma^2$. If the
functions involved are sufficiently smooth, one easily sees that the
integral is of order $\sigma^3$. If $r(t,s)$ vanishes in a quadratic
minimum in $s=t$ or elsewhere, then the integral is at most of order
$\sigma$. 

It is well known that the most probable path reaching $z$ at time $t$ is
represented by a straight line in the $(v,z)$-plane. Thus $r(t,s)$ vanishes
for some $s\neq t$ if and only if the most probable path reaching $d(t)$ as
already reached $d(s)$. In that case, there exists a time $u\in(0,t)$ such
that the tangent to the curve  $(v(s),d(s))_{s\geqs0}$ at $(v(u),d(u))$
goes through the origin, i.\,e.,  $(d(u)/v(u))'=0$. This situation can be
excluded under a convexity assumption on $d(t)$, which is equivalent to
Hypothesis H3. 
\end{remark}

The following lemma establishes the existence and some properties of a fixed
point of~\eqref{app3} under hypotheses tailored to our situation. We
will employ it in Section~\ref{ssec_fps} with 
$v(t) = \int_0^t \e^{2\alpha(s)}g(s)^2\6s$ and
$d(t) = (2-\deltap)\e^{\alpha(s)}$.

\begin{lemma}
\label{lem_appendix}
Assume that there are constants $\Delta, M_1,M_2>0$ such that the conditions 
\begin{align}
\label{app5a}
d(s)\vprime(s)-v(s)\dprime(s) &\geqs 
\Delta\vprime(s)\bigpar{1+\sqrt{v(s)}},\\
\label{app5b}
\abs{\tilde c(t,s)} &\leqs M_1, \\
\label{app5c}
M_2\vprime(s) &\geqs 1+v(s),
\end{align}
hold for all $0\leqs s\leqs t$. Then~\eqref{app2c} holds with a prefactor
$c(t)$ satisfying 
\begin{equation}
\label{app5d}
\abs{c(t) - c_0(t)} \leqs 
\frac{\eps}{1-\eps}
\frac{1+v(t)^{3/2}}{v(t)^{3/2}} \sup_{0\leqs s\leqs t}
\biggabs{\frac{v(s)^{3/2}}{1+v(s)^{3/2}}c_0(s)}, 
\end{equation}
whenever 
\begin{equation}
\label{app5e}
\eps := 2M_1 \biggpar{\frac{\e^{-\Delta^2/4\sigma^2}}\sigma \,t +
\frac{4M_2\sigma}{\Delta^2}} < 1. 
\end{equation}
\end{lemma}
\begin{proof}
We shall prove that $\cT$ is a contraction on the Banach space
$\cX$ of continuous functions $c:[0,t]\to[0,\infty)$, equipped with the
norm
\begin{equation}
\label{app6a}
\norm{c} = \sup_{0\leqs s\leqs
t}\biggabs{\frac{v(s)^{3/2}}{1+v(s)^{3/2}}c(s)}. 
\end{equation}
For any two functions $c_1, c_2\in\cX$, we have
by~\eqref{app5b} 
\begin{equation}
\label{app6}
\abs{\cT c_2(t) - \cT c_1(t)} \leqs  
\norm{c_2 - c_1} \frac{M_1}\sigma \int_0^t 
\frac{1+v(s)^{3/2}}{v(s)^{3/2}}\e^{-r(t,s)/2\sigma^2}\6s. 
\end{equation}
Using Assumption~\eqref{app5a}, we obtain
\begin{align}
\nonumber
\frac{d(s)}{v(s)} - \frac{d(t)}{v(t)} 
&= \int_s^t \frac{d(u)\vprime(u)-v(u)\dprime(u)}{v(u)^2}\6u\\
&\geqs \Delta \biggbrak{\frac{v(t,s)}{v(t)v(s)} +
2\frac{\sqrt{v(t)}-\sqrt{v(s)}}{\sqrt{v(t)v(s)}}},
\label{app7}
\end{align}
and thus 
\begin{equation}
\label{app8}
r(t,s) \geqs \Delta^2 \biggbrak{\frac{v(t,s)}{v(t)v(s)} +
4\frac{(\sqrt{v(t)}-\sqrt{v(s)})^2}{v(t,s)}} 
\geqs \Delta^2 \frac{v(t,s)}{v(t)} \biggbrak{1+\frac1{v(s)}}. 
\end{equation}
For the sake of brevity, we restrict our attention to the case
$v(t)>2$. We split the integral in~\eqref{app6} at times $s_1$ and
$s_2$ defined by $v(s_1)=1$ and $v(s_2)=v(t)/2$. By~\eqref{app5c}, the
integral on the first interval is bounded by 
\begin{equation}
\label{app8b}
2\int_0^{s_1} \frac1{v(s)^{3/2}} \e^{-\Delta^2/4\sigma^2v(s)}\6s
\leqs \frac{4M_2\sigma}{\Delta} \int_{\Delta^2/4\sigma^2}^\infty
\frac{\e^{-y}}{\sqrt{y}}\6y 
\leqs \frac{8M_2\sigma^2}{\Delta^2} \e^{-\Delta^2/4\sigma^2}. 
\end{equation}
The second part of the integral is smaller than $2t\e^{-\Delta^2/4\sigma^2}$
because $r(t,s)\geqs\Delta^2/2$ for $s_1<s<s_2$, while the last part
is bounded by  
\begin{equation}
\label{app9}
\int_{s_2}^t \e^{-\Delta^2 v(t,s)/2\sigma^2v(t)} \6s 
\leqs \frac{2^{3/2}M_2}{v(t)} \int_{s_2}^t \vprime(s) \e^{-\Delta^2
v(t,s)/2\sigma^2v(t)} \6s  
\leqs\frac{2^{5/2}M_2\sigma^2}{\Delta^2}. 
\end{equation}
This shows that $\cT$ is a contraction with contraction constant $\eps$,
and the result follows by bounding $\norm{c-c_0}=\norm{\cT^nc-\cT0}$ by a
geometric series. Here $0$ denotes the function which is zero everywhere.
\end{proof}

\begin{cor}
\label{cor_appendix}
Let the assumptions of Lemma~\ref{lem_appendix} be satisfied and
assume in addition that there exists a constant $M_3>0$ such that
\begin{equation}
\label{app9.1}
d(t)\vprime(t) - v(t)\dprime(t) \leqs M_3 (1+v(t)^{3/2})
\qquad\text{for all $t\geqs 0$.}
\end{equation}
Then
\begin{equation}
\label{app10}
c_0(t) \biggbrak{1 - \frac{\eps}{1-\eps}
\frac{M_2M_3}{\Delta}} \leqs 
c(t) \leqs c_0(t) \biggbrak{1 + \frac{\eps}{1-\eps}
\frac{M_2M_3}{\Delta}}
\end{equation}
holds for all $t>0$ such that $\eps = \eps(t)<1$, where $\eps$ is
defined by~\eqref{app5e}.
\end{cor}

\begin{proof}
The proof follows directly from the bounds~\eqref{app9.1} and
\begin{equation}
\label{app11}
\frac{1+v(t)^{3/2}}{v(t)^{3/2}}\frac{1}{c_0(t)} 
\leqs \frac{\sqrt{2\pi} M_2}{\Delta},
\end{equation}
the latter being a consequence of~\eqref{app5a} and~\eqref{app5c}.
\end{proof}


\small
\bibliography{../BG}
\bibliographystyle{amsalpha}               

\bigskip\bigskip\noindent
{\small 
Nils Berglund \\ 
{\sc FRUMAM, CPT--CNRS Luminy} \\
Case 907, 13288~Marseille Cedex 9, France \\
{\it and} \\
{\sc PHYMAT, Universit\'e de Toulon} \\
{\it E-mail address: }{\tt berglund@cpt.univ-mrs.fr}

\bigskip\noindent
Barbara Gentz \\ 
{\sc Weierstra\ss\ Institute for Applied Analysis and Stochastics} \\
Mohrenstra{\ss}e~39, 10117~Berlin, Germany \\
{\it E-mail address: }{\tt gentz@wias-berlin.de}
}


\end{document}